\documentclass[12pt]{conm-p-l}

\usepackage{tikz,enumerate,dsfont,amssymb,fullpage,wrapfig}
\usepackage{mathrsfs,mathtools, tipa}
\usetikzlibrary{decorations.pathmorphing}
\usetikzlibrary{decorations.markings}
\usetikzlibrary{arrows.meta,bending}
\usetikzlibrary{intersections}
\usetikzlibrary{arrows}
\usetikzlibrary{positioning}
\usetikzlibrary{cd}
\usetikzlibrary{fit}
\usepackage [all]{xy}
\SelectTips{cm}{12}
\usepackage{stmaryrd}
\usepackage{wasysym}
\usepackage{hyperref}
\usepackage[bbgreekl]{mathbbol}

\newtheorem{Thm}{Theorem}[section]
\newtheorem{Lem}[Thm]{Lemma}

\newtheorem {Cor}[Thm] {Corollary}

\theoremstyle{definition}
\newtheorem{Def}[Thm]{Definition}
\newtheorem{Exa}[Thm]{Example}

\theoremstyle{remark}

\numberwithin{equation}{section}

\DeclareMathSymbol\bbDelta  \mathord{bbold}{"01}

\newcommand{\amoon} {{\alpha_{\leftmoon}}}

\newcommand{\coker} {\operatorname{coker}}

\newcommand{\cyclicgroup}[1] {{C_{#1}}}

\newcommand{\teg}{{\{\pm 1\}}}
\newcommand{\TT}{T\!\rtimes\!\teg}
\newcommand{\mTT}{\mT\!\rtimes\!\teg}

\newcommand{\ul}[1]{\underline{#1}}

\newcommand{\ttr}{\mathbb Tr}

\newcommand{\bbL}{\mathbb L}

\colorlet{darkgreen}{green!60!black}
\colorlet{turquoise}{green!50!blue}
\colorlet{hellrot}{orange!50!red}
\colorlet{wein}{orange!30!red!80!blue}
\colorlet{ziegel}{hellrot!50!brown}

\newcommand {\specialmap} [4] {\text {$ #1\negmedspace : #2 #3 #4 $}}
\newcommand {\map} [3] {\specialmap {#1} {#2}{\to} {#3}}
\newcommand {\longmap} [3] {#1\negmedspace :#2 \longrightarrow #3}

\newcommand {\Aut}{\operatorname{Aut}}
\newcommand {\Hom} {\operatorname {Hom}}
\newcommand {\id} {\operatorname{id}}

\newcommand {\at}[1] {\arrowvert_{#1}}

\renewcommand {\(} {\left(}
\renewcommand {\)} {\right)} 
\newcommand {\spec} {\operatorname{spec}}

\newcommand {\CC} {\mathbb C}
\tikzstyle bigarrow=[postaction={decorate,decoration={markings,
    mark=at position 1 with {\arrow[scale=2]{>}}}}]

\newcommand{\mA}{{\mathcal A}}

\newcommand{\mC}{{\mathcal C}}
\newcommand{\mD}{{\mathcal D}}
\newcommand{\mE}{{\mathcal E}}

\newcommand{\mG}{{\mathcal G}}

\newcommand{\mI}{{\mathcal I}}

\newcommand{\mK}{{\mathcal K}}
\newcommand{\mL}{{\mathcal L}}
\newcommand{\mM}{{\mathcal M}}
\newcommand{\mN}{{\mathcal N}}
\newcommand{\mO}{{\mathcal O}}
\newcommand{\mP}{{\mathcal P}}
\newcommand{\mQ}{{\mathcal Q}}

\newcommand{\mS}{{\mathcal S}}
\newcommand{\mT}{{\mathcal T}}

\newcommand{\mV}{{\mathcal V}}

\newcommand{\mZ}{{\mathcal Z}}

\newcommand{\fh}{{\mathfrak h}}

\newcommand{\ft}{{\mathfrak t}}

\newcommand {\bbB} {{\mathbb B}}
\newcommand {\bbC} {{\mathbb C}}
\newcommand {\bbF} {{\mathbb F}}

\newcommand {\bbR} {{\mathbb R}}

\newcommand {\bbZ} {{\mathbb Z}}

\newcommand {\Bicat} {{\mathbf{Bicat}}}

\newcommand{\inv}{^{-1}}

\newcommand {\LE} {{\mathbf L}}
\newcommand {\Lv} {{\Lambda^{\!\vee}}}

\newcommand {\tensor}{\otimes}

\title{Looking for a Refined Monster}

\author{Nora Ganter}
\address{School of Mathematics and Statistics\\
  The University of Melbourne\\ Parkville VIC 3010\\ Wurundjeri Country\\ Australia}
\email{nganter@unimelb.edu.au}

\thanks{The author was supported by ARC grants
  DP160104912 and DP210103081.}


\date{\today} 


\begin{document}

\maketitle
\begin{abstract}
  We discuss some categorical aspects of the objects that appear in the construction of the Monster and other sporadic simple groups.
  We define the basic representation of the categorical torus $\mT$ classified by an even symmetric bilinear form $I$ and of the semi-direct product of $\mT$ with its canonical involution. We compute the centraliser of the basic representation  of $\mTT$ and find it to
  be a categorical extension of the extraspecial $2$-group with commutator $I\mod 2$. We study the inertia groupoid of a categorical torus and find that it is given by
  the torsor of the topological Looijenga line bundle, so that $2$-class functions on $\mT$ are canonically theta-functions.
  We discuss how discontinuity of the
  categorical character in our formalism means that the character of the basic representation fails to be a categorical class function. 
  We compute the automorphisms of $\mT$ and of $\mTT$ and relate these to the
  Conway groups.
\end{abstract}

\section{Introduction}
The largest of the sporadic finite simple groups in the Atlas
classification \cite{ATLAS85} is known by the name {\em the Monster}. Its
irreducible representations have dimensions
$$
  1,\quad 196883, \quad 21296876, \quad 842609326,\quad\dots
$$
These characters were found before the Monster was first constructed, and it was then that John McKay observed that the coefficients of the $j$-function
$$
  j(q)-744 = q^{-1}+196884q+ 21493760q^2+864229970q^3+\dots
$$
appeared to be counting the graded dimension of an infinite series of Monster
representations. The conjectural picture describing this phenomenon, which would
come to be known as Monstrous Moonshine, was further developed by Conway and
Norton \cite{ConwayNorton79} and generalised to involve pairs of commuting elements of the Monster
by Norton \cite{Mason87}. The conjecture was proved in the classical setting by Borcherds
\cite{Borcherds92} and in the generalised setting by Carnahan \cite{Carnahan07}, \cite{Carnahan10}, \cite{Carnahan12}.
The Monster itself was constructed by Griess
\cite{Griess82}, see also \cite{Tits85}. Its most conceptual construction
to date is due to Frenkel, Lepowsky and Meurman, who brought together the young disciplines of
infinite dimensional Lie theory and string theory, in order
to conceptualise the appearance of theta functions and integrate them into the construction of the Monster. This lead them to
develop the notion of a of vertex operator algebra, realising the Monster as the symmetry group of the Moonshine Module \cite{FrenkelLepowskyMeurman88}.
For a detailed account of the early history of
Monstrous Moonshine, we refer the reader to the introduction of
\cite{FrenkelLepowskyMeurman88}.

From the first encounters with the Monster, the subtle role of the number 24
was apparent. In the context of classical Moonshine, this manifests itself in the
level of the modular functions involved; in Generalised Moonshine, it manifests itself in the
occurrence of $24^{th}$ roots of unity and in the role of certain central extensions
of the centraliser subgroups of the Monster. It was conjectured by Mason
\cite{Mason02} that these were governed by a 3-cocycle of order $24$ in the group
cohomology of the Monster; this three-cocycle is called
the moonshine anomaly. We will write $M$ for the Monster group. In \cite{Johnson-Freyd19}, Johnson-Freyd studies a cocycle
$$\amoon\in H^3_{gp}(M;U(1)),$$
conjecturally that of Mason\footnote{Johnson-Freyd's moonshine anomaly is obtained using conformal nets, while Mason's 
  conjecture lives in the world of vertex operator algebras. The referee pointed out that, while these two formalisms are
  conjecturally equivalent to each other, the full VOA-conformal nets dictonary has yet to be worked out.}, and shows it to be of order $24$. 

In the days of the ATLAS project, Schur's work on the role of cocycles of
degree $2$ \cite{Schur04} \cite{Schur07} as defining central extensions was
already part of the standard repertoire of the group theorist.
The analogous classification of categorical central extensions by degree $3$
cocycles, on the other hand, was only just emerging, having been the topic of the unpublished
PhD thesis of Sính \cite{Sinh75}. It would be some time before
her {\em gr-cat\'egories}, nowadays better known under the name categorical
groups or, among non-group theorists, as $2$-groups\footnote{In this paper,
  the term $2$-group will be reserved for a finite group whose order is a power
  of $2$.}, became more mainstream. In any case, the ATLAS project was
concerned with the classification of finite simple groups; the community had no
interest in categorical groups. It is therefore not surprising that some of the
intrinsic categorical features of key objects in the construction of the sporadic
groups went unnoticed and, despite
their elementary nature, were not studied very deeply at the time.
In recent years, the question has been posed whether the Monster might, in fact,
occur in nature as a categorical group. Put differently, we know $\amoon$ to
classify a categorical group $\mM$ sitting inside an exact sequence
$$
\begin{tikzcd}
  1\ar[r]&\bbB \cyclicgroup{24} \ar[r] & \mM\ar[r] &M\ar[r] & 1,
\end{tikzcd}
$$
where $\cyclicgroup{24}$ denotes the cyclic group of order $24$, 
and $\bbB \cyclicgroup{24}$ is the one object category with automorphism group $\cyclicgroup{24}$.
We can see what appears to be a shadow of $\amoon$ in the Monster's natural habitat: Moonshine
and the Moonshine module. Defining a categorical group from a $3$-cocycle,
however, is the most tedious and unintuitive construction one could imagine,
involving triples of elements of the Monster group in the case of $\mM$.
One would instead hope to realise $\mM$ as the automorphisms of an object of a
bicategory. To be interesting, this object should be natural and easily
constructed and have the potential to simplify our understanding of the Monster.
The paper at hand will not provide an answer to this problem. Instead, I will
try to give an idea of where to look and how to think about these structures.

\medskip
Let $(\Lambda,I)$ be an even unimodular lattice, and let $O(\Lambda,I)$ be its
group of isometries. These data are already of a categorical nature: let
$\ft = \Lambda\tensor_\bbZ \bbR$ 
  and consider the torus 
$$ T \,=\, \Lambda\tensor_\bbZ U(1)  \,\cong\, \ft\,/\,\Lambda, 
$$
where on the right $\Lambda$ is identified with $\Lambda\tensor1$.
Then the Lie group cohomology \cite{WagemannWockel15} of $T$ is computed by the isomorphism
$$
  H^3_{gp}(T;U(1))\cong H^4(BT;\bbZ)
$$
\cite[Corollary 97]{Schommer-Pries11}. The elements of 
$H^4(BT;\bbZ)^{O(\Lambda,I)}$
are exactly  the even symmetric bilinear forms on $\Lambda$ invariant under
the action of $O(\Lambda,I)$, see
for instance \cite{Ganter18}. In particular, $I$ itself defines such an element.
As in the finite case, elements of the third Lie
group cohomology classify categorical central extensions.
The categorical
torus
$$
\begin{tikzcd}
  1\ar[r]&\bbB U(1) \ar[r] & \mT\ar[r] &T\ar[r] & 1
\end{tikzcd}
$$
corresponding to $I$ under this classification is easily constructed. We review two equivalent constructions in
Sections \ref{sec:categorical_tori} and \ref{sec:multiplicative_gerbe} below.
The motivating example is the categorical Leech torus $\mT_{Leech}$, where $(\Lambda,I)$ is taken to be the Leech lattice with its usual biliear form. In this setting, the isometry group is Conway's group $O(\Lambda,I)=Co_0$.

In the paper at hand, we shall begin to explore the representation and character theory of categorical tori. While far from presenting a mature theory, we will see glimpses of how the categorical picture captures important features, such as a theta function formalism, from a geometric point of view. At the same time, it is
analytically much less intricate than loop group representations, which  have long been studied as geometric counterpart to representations of affine Lie algebras,
such as, for instance, the (untwisted) affine Heisenberg algebra
$\widehat\fh$ in \cite{FrenkelLepowskyMeurman88}.
The emerging philosophy views categorical Lie groups as a third pillar of the theory. The relationship between loop groups and categorical groups has been studied, for instance, in \cite{BaezStevensonCransSchreiber07}, \cite{Waldorf12}, \cite{BrylinskiMcLaughlin94}, \cite{LudewigWaldorf23}, \cite{Ganter18}, and there are already instances, where the finite dimensional categorical picture sheds light on the choices one ought to make in the loop group setting.

We will be interested in symmetries. The invariance of the cohomology class $I$ under the action of its isometry group $O(\Lambda,I)$ realises the centre preserving automorphisms of the categorical torus $\mT$ classified by $I$ as a categorical central extension of $O(\Lambda,I)$. This computation can also be found in \cite{Waldorf22}.
We will also be interested in the categorical group
$\mT\rtimes\{\pm1\}$, which can be thought of as a semi-direct product and plays the role of `orbifoldisation'.
The outer automorphisms of $\mathcal T\rtimes\{\pm1\}$ form a categorical extension of \[PO(\Lambda,I)=O(\Lambda,I)/\{\pm\id\}.\] In the case of the Leech lattice, this group is the Conway group $Co_1$. A categorical central extension of the extraspecial $2$-group $2^{1+24}$ in ATLAS notation turns up as the centraliser of the basic representation of $\mT_\bbC\rtimes\{\pm 1\}$. It is therefore expected that the normaliser of this basic representation is closely related to the subgroup $2^{1+24}.Co_1$ of the Monster.
\subsection{Acknowledgments}
It is a pleasure to thank Matthew Ando and Gerd Laures, John McKay and Igor Frenkel, Geoffrey Mason, Konrad Waldorf and Thomas Nikolaus, Theo Johnson-Freyd, and David Treumann for helpful conversations. Special thanks go to Igor Frenkel for convincing me to bring this circle of ideas to paper.
Many thanks go to the referee and the referee's PhD student for extensive feedback and help in reorganising the exposition. 
\section{Background on crossed modules and strict categorical groups}\label{sec:background}
\subsection{Crossed modules and strict categorical groups}\label{sec:background}
By a {\em categorical group} we will mean a small monoidal groupoid with weakly invertible objects. A categorical group is called {\em strict} if its objects are
strictly invertible and all structure morphisms (units and associators) are identity maps. Strict categorical groups are most easily described using crossed
modules.
We will largely follow the
conventions in \cite{Noohi07}.
\begin{Def} \label{def:crossed}
  A crossed module $\Psi$ consists of two groups $G_0$ and $G_1$, equipped
  with an action of $G_0$ on $G_1$ from the right, denoted
  \begin{eqnarray*}
      G_1\times G_0 &\longrightarrow& G_1\\
  	  (\alpha,x) &\longmapsto &\alpha^x
  \end{eqnarray*}
  and a group homomorphism $\map\psi{G_1}{G_0}$, 
 written $\alpha\mapsto \underline\alpha$, such that the following two
  conditions are satisfied:
  \begin{enumerate}
    \item[(CM1)] \quad $\underline{\alpha^x} = x\inv\underline\alpha
      x$\quad (equivariance)
      \smallskip
      
    \item[(CM2)] \quad $\beta^{\underline{\alpha}} =
      \alpha\inv\beta\alpha$\quad (Pfeiffer identity).
  \end{enumerate}
\end{Def}
Associated to such a crossed module $\Psi$ is the four term exact sequence
\begin{equation*}
  \label{eq:crossed_extension}
  \begin{tikzcd}
    0\ar[r] & \ker(\psi)\ar[r] & {G_1}\ar[r,"\psi"] & {G_0}\ar[r] & \coker(\psi)\ar[r] & 0,
  \end{tikzcd}
\end{equation*}
called the {\em crossed module extension} corresponding to $\Psi$.
\begin{Def}\label{def:homomorphism}
  A homomorphism from a crossed module $\Psi$ to another crossed module $\Phi$ is a commuting square
  \[
    \begin{tikzcd}
      G_1\ar[r,"f_1"] \ar[d,swap,"\psi"] & H_1\ar[d,"\phi"]\\
      G_0\ar[r,"f_0"] & H_0,
    \end{tikzcd}
  \]
  where the vertical maps are the structure maps of the two crossed modules $\Psi$ and $\Phi$, and the horizontal maps
  $f_0$ and $f_1$ are group homomorphisms that are compatible with the actions in the sense that $$f_1(\alpha^x) = f_1(\alpha)^{f_0(x)}.$$ 
\end{Def}
\begin{Def}\label{def:S(psi)}
In the situation of the Definition \ref{def:crossed}, the strict categorical group associated to 
$\Psi$ is the groupoid 
\begin{displaymath}
    \begin{tikzpicture}
       \node at (-1.8,-.5)    {$\mS(\Psi)\,\,\,\,\simeq\,$};
    	\node at (0,0.2) {$G_0\ltimes G_1$};
    	\node at (0,-1.2) {$G_0$};
    	\draw[->] (-.25,0) -- (-.25,-1);
    	\draw[->] (.25,0) -- (.25,-1);
    \end{tikzpicture}
\end{displaymath}
with source $pr_1$ and target $pr_1\cdot \psi$,
equipped with the strict monoidal structure given by the group multiplications.
\end{Def}
The construction of Definition \ref{def:S(psi)} defines a functor $\mS$ from the category of crossed modules and crossed module homomorphisms to the category
of strict categorical groups and strict monoidal functors.
If $\mS$ is a categorical group, and $1$ is its unit object, then we will write $\pi_0(\mS)$ for the group of isomorphism classes of $\mS$ and $\pi_1(\mS) = Aut_\mS(1)$
for the Bernstein centre of $\mS$.
One has $\pi_1(\mS(\Psi))=\ker(\psi)$ and
$\pi_0(\mS(\Psi)) = \coker(\psi)$.
When we speak of a categorical group extension of the form
\[
  \begin{tikzcd}
    1\ar[r] & \bbB A\ar[r]& \mG\ar[r] & G\ar[r] & 1,
  \end{tikzcd}
\]
where $G$ is a group and $A$ is an abelian group, we mean a categorical group $\mG$ with $\pi_0(\mG)\cong G$ and $\pi_1(\mG)\cong A$.
\subsection{Strict automorphisms and centre of a crossed module}
Let $\Psi$ be a crossed module as in Definition \ref{def:crossed}.
Equip the set of crossed homomorphisms
  \[
    Cross(G_0,G_1) = \{\map\chi{G_0}{G_1}\mid \chi(xy)\,\,=\,\, \chi(x)^y\chi(y)\}
  \]
  with the composition rule
  \begin{equation}
    \label{eq:vertical_composition}
    \(\chi_1\circ\chi_2\)(x) = \chi_2(x)\cdot\chi_1(x\cdot\ul{\chi_2(x)}),
  \end{equation}
  and write $Cross^\times(G_0.G_1)$ for the group of invertible elements in this semi-group.
  In the Lie setting, we will require crossed homomorphisms to be analytic.
\begin{Def}[{\cite{Norrie90} \cite{Breen92}}]
\label{def:Norrie_actor}
  The {\em actor crossed module} $Act(\Psi)$ consists of the groups $Act_0(\Psi) = Aut(\Psi)$ of crossed module automorphisms and 
  $Act_1(\Psi) = Cross^\times(G_0,G_1)$, along with the group homomorphism
  $\delta\negmedspace: Act_1\to Act_0$ sending the crossed homomorphism $\chi$ to the pair of group automorphisms $x\mapsto x\ul{\chi(x)}$ and $\alpha\mapsto\alpha\chi(\ul\alpha)$. 
  The action of $Act_0$ on $Act_1$ is given by $$\chi^{(f_0,f_1)} = f_1\inv\chi f_0.$$ By the adjoint representation of $\Psi$ we will mean the
  crossed square 
\[
  \begin{tikzcd} 
    & \alpha\ar[rr,mapsto] &[-2ex]& \chi_\alpha & \\[-3ex]
    & {G_1}\ar[dd,swap,"\psi"] \ar[rr] && {Cross^\times(G_0,G_1)}\ar[dd,"\delta"]  &
     \\[-2ex]  \boldsymbol{Ad:}  \\[-3ex]
    & {G_0}  \ar[rr] && Aut(\Psi) 
    \\[-3ex]
    & x\ar[rr,mapsto] && \(x(-)x\inv,(-)^{x\inv}\),
  \end{tikzcd}
\]
where $\chi_\alpha(x) =  \alpha^x\alpha\inv$. We will also write $c_x=Ad_0(x)$ for conjugation with $x$.
\end{Def}
We will mostly be working with a variation of Definition \ref{def:Norrie_actor} that considers only the centre-preserving crossed module homomorphisms; i.e., 
we will work with the crossed module $Act^+(\Psi)$ and the crossed square $Ad^+$ obtained by replacing $Aut(\Psi)$ with the automorphisms of $\Psi$ that act as the identity on $\pi_1(\Psi)$.

To define the Drinfeld centre of a crossed module, we follow Pirashvili\footnote{Note that our crossed module conventions differ from those in the work of Norrie and Pirashvili referenced here, who both use left actions. A variation of the definition of centraliser
  also appears in \cite[Page 41]{CarrascoGarzon04}.}, see \cite{pirashvili2023centre}.
Given $\Psi$ as in Definition \ref{def:crossed}, we form the semi-direct product
$$G_0\ltimes Cross(G_0,G_1),$$
where the action of $G_0$ is via $Ad_0$ and the action in $Act(\Psi)$, so, $$\chi^x(z) = \chi(x z x\inv)^x.$$
Consider the subgroup
$$Z_0(\Psi) = \left\{(x,\xi)\in G_0\ltimes Cross(G_0,G_1) \mid Ad_0(x\inv) = \delta(\xi)\right\}.$$
In $Z_0(\Psi)$, the multiplication in the semi-direct product simplifies to
\[(x,\xi)(y,\upsilon)  = \(xy,z\mapsto \upsilon(z)\xi(z)^y\),\]
and in particular, $Z_0(\Psi)$ is a group.
\begin{Def}[{\cite{pirashvili2023centre}}]
  The centre of $\Psi$ is the crossed module $Z(\Psi)$ given by the map
  \begin{eqnarray*}
    \zeta\negmedspace : G_1 & \longrightarrow & Z_0(\Psi)\\
    \alpha & \longmapsto & (\ul\alpha,Ad_1(\alpha\inv)),
  \end{eqnarray*}
  together with the action $\alpha^{(x,\xi)} = \alpha^x$.
\end{Def}
The crossed module $Z(\Psi)$ encodes the Drinfeld centre of $\mS(\Psi)$,
\[
  \mS(Z(\Psi)) \cong \mZ(\mS(\Psi)).
\]
\subsection{Weak maps between crossed modules}
\label{app:Conventions}
In \cite{Noohi07}, Noohi worked out how (not necessarily strict) monoidal functors and monoidal natural transformations between strict categorical groups are encoded in terms of the relevant crossed modules. All of Noohi's crossed modules are discrete. In the Lie setting, the butterfly formalism of \cite{AldrovandiNoohi09} is the more suitable notion of weak morphism. Since all the relevant groups in our examples are homotopically discrete, the two formalisms turn out to be equivalent; we do not concern ourselves with the distinctions here and refer the reader to \cite{Waldorf22} for a detailed comparison of different notions in the Lie case. 
\begin{Def}[{compare \cite[Definition 8.4]{Noohi07}}]
  \label{def:weak_morphism}
	Let $\Phi$ and $\Psi$ be crossed modules with structure maps
  $\longmap\phi{H_1}{H_0}$ and $\longmap\psi{G_1}{G_0}$. A weak morphism from $\Phi$ to $\Psi$ consist of a pair of pointed set maps $\longmap{p_i}{H_i}{G_i}$, for
	$i\in\{0,1\}$, together with a map
	\[\longmap\kappa{H_0\times H_0}{G_1}\]
	satisfying
  \[
  \begin{array}{lccrcll}    
    (W1)&\quad&\quad&
    p_0(\ul\alpha)&\,=\,&\ul{p_1(\alpha)}& \\[+8pt]

    (W2)&&&
    p_0(x)p_0(y)\ul{\kappa_{x,y}}&\,=\,&p_0(xy) &\quad\quad\quad
                                         \text{(K\"unneth 1)}\\[+8pt]
    (W3)&&&
    p_1(\alpha)p_1(\beta)\kappa_{\ul\alpha,\ul\beta}&\,=\,&p_1(\alpha\beta)
                                               &             \quad\quad\quad \text{(K\"unneth 2)} \\[+8pt]
    (W4)&&&\kappa_{x,y}^{p_0(z)}\kappa_{xy,z}&\,=\,&
    \kappa_{y,z}\kappa_{x,yz}&\quad\quad\quad\text{(cocycle), and}\\[+8pt]
    (W5)&&&
    (p_1(\alpha))^{p_0(x)}\kappa_{\ul\alpha,x}&\,=\,&p_1(\alpha^x)\kappa_{x,\ul{\alpha^x}}
                        &\quad\quad\quad\text{(equivariance).}   
  \end{array}
  \]
  We will use the notation $\kappa_{x,y,z}$  to refer to either side of (W4).
\end{Def}
Combined with the other four, our condition (W5) is equivalent to that
in \cite{Noohi07}. Note that (W3) implies $\kappa_{1,1}=1$, and hence
(W4) implies $\kappa_{x,1}=1=\kappa_{1,y}$. It follows that the restriction of
$p_1$ to $\pi_1(\mS(\Psi))$ is a group homomorphism, and so is
$\pi_0(p_0)$.
Horizontal composition of weak morphisms is given by
\begin{equation}
  \label{eq:composed_kuenneth}
  (f,\kappa)(g,\gamma) = (f\circ g,\beta)
\end{equation}
with
\[
  \beta_{x,y} \,=\, f_1(\gamma_{x,y})\kappa_{g_0(x),g_0(y), \ul{\gamma_{x,y}}}.
\]
It is a tedious yet straightforward computation to check that $\beta$ satisfies the cocycle condition (W4), and we have
\[
  \beta_{x,y,z} \,=\, f_1(\gamma_{x,y,z})\kappa_{g_0(x),g_0(y),g_0(z), \ul{\gamma_{x,y,z}}}.
\]
Associativity and well-definedness of \eqref{eq:composed_kuenneth} follow.
  The monoidal functor corresponding to $(p,\kappa)\in wHom(\Phi,\Psi)$ is the functor
  \[
    \longmap P{\mS(\Phi)}{\mS(\Psi)}
  \]
  which is given by $p_0$ on objects, by
  \[
    (x,\alpha)\longmapsto(p_0(x),p_1(\alpha)\kappa_{x,\ul\alpha})
  \]
  on arrows and has K\"unneth isomorphisms 
  \[
    K_{x,y} \,=\, (p_0(x)p_0(y),\kappa_{x,y}).
  \]
\begin{Def}
  Given a crossed module $\Psi$ as in Definition \ref{def:crossed}, we define the weak actor crossed module ${wAct(\Psi)}$ of $\Psi$ as follows:
  \[
    wAct(\Psi)_0 \,=\, \{(f,\kappa)\in wHom(\Psi,\Psi) \mid \text{$f_0$ and $f_1$ are bijections}\}
  \]
  is the group of invertible weak endomorphisms of $\Psi$ with the multiplication \eqref{eq:composed_kuenneth},
  while
  \[
    wAct(\Psi)_1 \,=\, Maps_*(G_0,G_1)^\times
  \]
  is the group invertible elements in the semi-group of pointed maps from $G_0$ to $G_1$ with multiplication given by \eqref{eq:vertical_composition}.
  The structure map $\longmap{\delta}{wAct(\Psi)_1}{wAct(\Psi)_0}$ sends $\eta$ to $(f,\kappa)$ with
  \begin{eqnarray*}
    f_0(x) & = & x\ul{\eta(x)}\\
    f_1(\alpha) & = & \alpha\eta(\ul{\alpha})\\
    \kappa_{x,y} & = & \eta(y)\inv\(\eta(x)^y\)\inv\eta(xy).
  \end{eqnarray*}
  The action of $wAct(\Psi)_0$ on $wAct(\Psi)_1$ is given by
  \begin{eqnarray*}
    \eta^{(f,\kappa)}(x) & = & f_1\inv\(\eta\(f_0(x)\)\)\, \gamma_{f_0(x),\ul{\eta(f_0(x))}}\\[+3pt]
     & = & f_1\inv\(\eta\(f_0(x)\)\kappa_{x,\ul{\eta^{(f,\kappa)}(x)}}\inv\),
  \end{eqnarray*}
  where in the first row, $(f\inv,\gamma) = (f,\kappa)\inv$, and in the second row,
  \[\ul{\eta^{(f,\kappa)}(x)} = x\inv f_0\inv\(f_0(x)\ul{\eta(f_0(x))}\).\]
\end{Def}
We will call a weak endomorphism of $\Psi$ {\em centre-preserving} if its restriction to $\pi_1(\Psi)$ is the identity.
Since all the weak automorphisms in the image of $\delta$ are centre preserving,
the centre-preserving weak automorphisms of $\Psi$ form a sub-crossed module of $wAct(\Psi)$, which we will denote $wAct^+(\Psi)$.
The strict categorical group $\mS(wAct(\Psi))$ encodes the strictly invertible weakly monoidal endofunctors of $\mS(\Psi)$.
We will will use the notation\footnote{Note the notation clash with \cite{CarrascoGarzon04}, where $\mA ut$ is reserved for strict automorphisms, while our $\mS(wAct(\Psi))$ is contained in their $\mE q(\mG)$.} $\mA ut(\mS(\Psi))$ for this categorigal group and $\mA ut^+(\mS(\Psi))$ for $\mS(wAct^+(\Psi))$.
\subsection{Automorphisms} 
Let $G$ be a group. Then the crossed
module extension
\begin{equation}\label{eq:Ad_groups}
\begin{tikzcd}
  1\ar[r]& Z(G) \ar[r]& G\ar[r,"Ad"] & Aut(G) \ar[r]& Out(G)\ar[r] & 1\\[-.65cm]
  &&g\ar[r,mapsto]&c_g=g(-)g\inv, & 
\end{tikzcd}
\end{equation}
encodes, via the construction of Definition \ref{def:S(psi)}, the categorical group
$$\mA ut(\bbB G)\simeq \mS(Ad).$$
\begin{Exa}
  For $G=\TT$ the extension \eqref{eq:Ad_groups} takes the form
  \[
    \begin{tikzcd}
      1\ar[r]&[-1ex] Z(\TT) \ar[r]\ar[d,equal]&[-1ex] \TT \ar[r,"Ad"]\ar[d,equal] &
      Aut(\TT) \ar[r]\ar[d,equal]&[-1ex] Out(\TT)\ar[r]\ar[d,equal] &[-1ex] 1\\
      1\ar[r]& T[2] \ar[r]& \TT\ar[r,"(-)^2"] &
      T\!\rtimes\! GL(\Lv) \ar[r]& PGL(\Lv)\ar[r] & 1,
    \end{tikzcd}
  \]
  where $PGL(\Lv)=GL(\Lv)/\{\pm id\}$.
\end{Exa}
Given a strict categorical group $\mG=\mS(\Psi)$, we consider the composite
\[
  \begin{tikzcd}
    \Psi\ar[r,"\boldsymbol{Ad^+}"] &[+1ex] Act^+(\Psi)\ar[r] & wAct^+(\Psi),
  \end{tikzcd}
\]
where $\boldsymbol{Ad}$ is the crossed square of Definition \ref{def:Norrie_actor}.
Applying $\mS$, we obtain the categorical
crossed module
\cite{Breen92}
\cite{CarrascoGarzonVitale06} 
\[
  \begin{tikzcd}
    \mA d^+\negmedspace : \mG \ar[r] &
    \mA ut^+(\mG),
  \end{tikzcd}    
\]
encoding the automorphism 3-group
\[
  \mA ut^+(\bbB\mG)\,\,=\,\, \mA ut^+(\mG) /\negmedspace/ \mG
\]
\cite[4]{Rousseau03}
with 2-isomorphism classes
\[
  \mC oker(\mA d^+) \,\, \simeq \,\, \mO ut^+(\mG)
\]
(see \cite[p. 42]{CarrascoGarzon04}), and Bernstein centre equal to the Drinfeld
centre of $\mG$,
\[
  \mK er(\mA d^+) \,\, \simeq \,\, \mZ(\mG).
\]
Its homotopy groups fit into the exact sequence
  \begin{equation}
    \label{eq:seven-terms}
    \hspace{-3.4cm}
    \begin{tikzcd}
      1\ar[r]& (\pi_1\mG)^{\pi_0\mG}\ar[r]&
      \pi_1\mG\ar[r]&
      Cross(\pi_0\mG,\pi_1\mG)  
      \arrow[dll,controls={+(7,-1) and +(-6,1)}]\\[+2ex]
      &\pi_0\mZ\mG\ar[r]&
      \pi_0\mG\ar[r]&
      \pi_0\mA ut^+(\mG)\ar[r]&
      \pi_0\mO ut^+(\mG) \ar[r]& 1
    \end{tikzcd}
\end{equation}
\cite[Corollary 3.3]{CarrascoGarzon04}.
It is sometimes helpful to restrict our attention to the automorphisms over $G=\pi_0(\mG)$, i.e., the automorphisms of the extension.
Given a central extension of finite dimensional Lie groups
\[
  \begin{tikzcd}
    1\ar[r] &A\ar[r] &\widetilde G\ar[r,"\pi"] & G\ar[r]& 1
  \end{tikzcd}
\]
defined by a locally continuous $A$-valued 2-cocycle $\beta$ on $G$, the automorphisms of $\pi$ fixing the central subgroup are identified with the first cohomology group
\[
  H^1_{gp}(G;A) \,=\, Hom(G,A).
\]
We therefore have a commuting diagram with exact rows
\begin{equation}
  \label{eq:Aut(pi)}
  \begin{tikzcd}
    1\ar[r] &Z(\widetilde{G})\ar[r]\ar[d] &\widetilde
    G\ar[r]\ar[d,"Ad"] & Inn(\widetilde G)\ar[r]\ar[d,hook]& 1 \\ 
    1\ar[r] &H^1_{gp}(G;A)\ar[r] & Aut^+(\widetilde G)\ar[r] & Stab([\beta])\ar[r]& 1,
  \end{tikzcd}
\end{equation}
whose last term refers to the stabilizer of the class $[\beta]$ in
$Aut(G)$. 
Let now $\mQ$ be the multiplicative $A$-bundle gerbe over $G$, corresponding to the categorical central extension
\[
  \begin{tikzcd}
    \bbB A\ar[r] & \mG\ar[r] & G,
  \end{tikzcd}
\]
where $\mG$ is a strict categorical group.
Then the 
auto-equivalences of $\mQ$ are identified with the auto-equivalences
of $\mG$ over $G$ that fix the centre.
In other words, the categorical group $\mA ut(\mQ)$
is the kernel of the homomorphism
\[
  \begin{tikzcd}
    Cl\negmedspace : \mA ut^+(\mG)\ar[r] & Aut(G) \\[-3.5ex]
    \phantom{Cl\quad\quad:}f\ar[r,mapsto] & \pi_0(f)
  \end{tikzcd}
\]
sending an auto-equivalence $F$ to its effect on isomorphism classes.
If $\mQ$ is classified by the \v Cech-simplicial 3-class
$[\alpha]$ then the image of $Cl$ consists of the stabilizer
of $[\alpha]$ inside $Aut(G)$,
\[
  im(Cl) \,=\, Stab([\alpha]).
\]
Further, $\mA ut(\mQ)$ is identified with the strict categorical group
associated to
\[
  \begin{tikzcd}
    \check C^1(G;A) \ar[r,"d"] & \check Z^2(G;A).
  \end{tikzcd}
\]
The universal property of kernel induces the dashed arrow in the
diagram
\[
  \begin{tikzcd}
    \mG\at{Z(G)}\ar[r]\ar[dd,dashed
    ] & \mG\ar[rr, two
  heads, "Cl\circ  \mA d"]\ar[dd,"\mA d^+"]&& 
    Inn(G)\ar[dd,hook] \\ \\
    \mA ut(\mQ)\ar[r] & \mA ut^+(\mG)\ar[rr,"Cl"] && Aut(G),
  \end{tikzcd}
\]
and we obtain isomorphisms
\[
  \pi_1(\mG\at{Z(G)}) \,\cong\, \pi_1(\mG)
\]
and
\[
  H^1_{gp}(G;A))\,\cong\,\pi_1(\mA ut^+(\mG))
\]
and a map of short exact sequences
\begin{equation}
  \label{eq:Z->G->Inn}
  \begin{tikzcd}
    1\ar[r] & Z(G) \ar[d]\ar[r] & G \ar[r,"Ad"]\ar[d,"\pi_0(\mA d^+)"] &
    Inn(G)\ar[d,hook] \ar[r] & 1\\ 
    0\ar[r]&H^2_{gp}(G,A) \ar[r]& \pi_0(\mA ut^+(\mG))\ar[r] &
    Stab([\alpha])\ar[r]&1. 
  \end{tikzcd}
\end{equation}
\section{The basic representation of a categorical torus}\label{sec:basic_representation}
\subsection{Categorical tori}
\label{sec:categorical_tori}
Let $I$ be an even symmetric
bilinear form on a lattice $\Lv$, not necessarily unimodular. We write $\Lambda=\Hom(\Lv,\bbZ)$ for the dual lattice of $\Lv$, and view $\Lambda$ as the weight lattice of the torus $T=\Lv\tensor_\bbZ U(1)$ with Lie algebra $\ft=\Lv\tensor_\bbZ \bbR$.
Let $J$ be any choice of integral bilinear form on $\Lambda$ satisfying \label{page:J}
$$I(m,n) = J(m,n)+J(n,m).$$
By mild abuse of notation, we write $J$ also for the
bilinear form $J_\bbR=J\tensor\bbR$ on $\ft$. 
Then we have the crossed module $\Theta$ given by
\begin{eqnarray*}
  \theta \negmedspace : \Lv\times \bbR/\bbZ & \longrightarrow & \ft\\
  (m,z) & \longmapsto & m,
\end{eqnarray*}
where $x\in\ft$ acts on $\Lv\times\bbR/\bbZ$ via
$$
  (m,[a])^x  =  (m,[a+J(m,x)])
$$
(c.f.\ \cite{Ganter18}). The categorical torus classified by $I$ is the strict categorical group $\mT = \mS(\Theta)$.
So, $\mathcal T$ has objects $\ft$ and arrows
$$
\begin{tikzpicture}
  \node at (0,0) (x) {$x$};
  \node at (2.5,0) (xm) {$x+m$,};
  \draw[->] (x) -- node [midway,above] {$[a]$} (xm);
\end{tikzpicture}
$$
with source $x\in\ft$, target $x+m$ with  $m\in\Lv$, and label $[a]\in\bbR/\bbZ$. Composition is given by addition of labels
$$
\begin{tikzpicture}
  \node at (0,0) (x) {$x$};
  \node at (2.5,0) (xm) {$x+m$,};
  \node at (6,0) (xmn) {$x+m+n$,};
  
  \draw[->] (x) -- node [midway,above] {$[a]$} (xm);
  \draw[->] (xm) -- node [midway,above] {$[b]$} (xmn);
  \draw[->,bend right=25] (x) to node [midway,below] {$[a+b]$} (xmn);  
\end{tikzpicture}
$$
and the monoidal structure is given by
$$
\begin{tikzpicture}
  \node at (0,0) (x) {($x$};
  \node at (2.5,0) (xm) {$x+m)$};

  \node (bullet) [right of=xm] {$\bullet$};
  \draw[->] (x) -- node [midway,above] {$[a]$} (xm);

  \node [right of=bullet,xshift=-.4cm] (y) {$(y$};
  \node at (6.5,0) (yn) {$y+n$)};
  \draw[->] (y) -- node [midway,above] {$[b]$} (yn);

  \node (eq) [right of=yn] {$=$};
  \draw[->] (x) -- node [midway,above] {$[a]$} (xm);

  \node [right of=eq] (xy) {$(x+y$};
  \node at (14.5,0) (xymn) {$x+y+m+n)$,};
  \draw[->] (xy) -- node [midway,above] {$[a+b+J(m,y)]$} (xymn);

\end{tikzpicture}
$$
with trivial associators and unit maps. We will sometimes find it convenient to identify the circle group $\bbR/\bbZ$ with $U(1)$ via the exponential map
$a\mapsto e^{2\pi i a}$ and to write the labels multiplicatively.

We similarly
construct $\mT\rtimes\{\pm1\}$ as the strict categorical group associated to the crossed module $\Theta'$ with structure map
\[\theta'\negmedspace :\Lv\times\bbR/\bbZ\longrightarrow \ft\rtimes\{\pm1\},\] which is obtained by replacing $\ft$ with $\ft\rtimes\{\pm1\}$ in the construction of $\Theta$.
  In this setting the action of $\mathfrak t$ on arrows is as before and $-1$ acts on
  everything by the involution sending 
  $x\in\ft$ to $-x$ and sending the
  arrow
  \begin{equation}
    \label{eq:iota}
    x\xrightarrow{\,\,\,\,[a]\,\,\,\,}x+m \text{\quad\quad\quad to \quad\quad\quad}
    -x\xrightarrow{\,\,\,\,[a]\,\,\,\,}-x-m.
\end{equation}
  Finally, we define the complexifications $\Theta_\bbC$ of $\Theta$ and $\Theta'_\bbC$ of $\Theta'$ by replacing $\ft$ with the complex Lie algebra
  $\fh=\Lv\tensor_\bbZ\bbC$ and replacing the circle group $U(1)$ with the group of complex units $\bbC^\times$. These crossed modules define the complex categorical
  torus $$\mT_\bbC=\mS(\Theta_\bbC)$$ and the semi-direct product $$\mT_\bbC\rtimes\{\pm1\}=\mS(\Theta'_\bbC).$$
\subsection{The basic representation of a categorical torus}
We write 
$$\widehat T = \Hom(T,U(1))\quad\text{ and }\quad \check T = \Hom(U(1),T)$$
for the character and cocharacter lattice of $T$. These are isomorphic to $\Lambda$ and $\Lv$, and their elements are written multiplicatively. We write $e^{2\pi i\lambda}\in\widehat T$ for the character corresponding to the weight $\lambda$.  
The basic representations of $\mT$ is most easily defined as a strict map of crossed modules. Let 
\[
  \bbC[\widehat T] = \bbC \{e^{2\pi i\lambda}\}_{\lambda\in\Lambda}
\]
\[
  e^{2\pi i\lambda}\cdot e^{2\pi i\mu} = e^{2\pi i(\lambda+\mu)}
\]
be the group algebra of $\widehat T$, and similarly for $\check T$.
Let $\Aut_{var}(T_\bbC)$ be the group of automorphisms of the complex variety
$$
  T_\bbC = \spec(\bbC[\widehat T]).
$$
Viewing $\bbC[\widehat T]$ as the algebra of polynomial functions on $T_\bbC$, the group $\Aut_{var}(T_\bbC)$ acts on $\bbC[\widehat T]$ by precomposition.  
Before we state the definition, we introduce some further notation: we will write
\begin{eqnarray*}
  J^\sharp\negmedspace :\Lv & \longrightarrow & \Lambda \\
    m & \longmapsto & J(m,-)
\end{eqnarray*}
and
\begin{eqnarray*}
J^\flat\negmedspace :  \Lv & \longrightarrow & \Lambda \\
   m  & \longmapsto & J(-,m)
\end{eqnarray*}
for the adjoints of $J$, and
$$\exp\negmedspace :\fh\longrightarrow T_\bbC$$ for the exponential map $$\exp=\Lv\tensor_\bbZ e^{2\pi i(-)}.$$
\begin{Def}\label{def:basic_rep}
  The basic representation of $\Theta'_\mathbb C$ is the strict crossed module homomorphism   
  \begin{center}
    \begin{tikzpicture}
      \node at (-1.75,-2) {$r_{bas}\negmedspace:$};
      \node at (0,0) [name=a] {$(m,z)$};
      \node at (3,0) [name=b,anchor=west] {$ze^{2\pi iJ^\sharp(m)}$};
      \node at (0,-1) [name=c] {$\Lv\times {\bbC^\times}$};
      \node at (3,-1) [name=d,anchor=west] {$\CC[\,\widehat T\,]^\times$};
      \node at (0,-3) [name=e] {$\fh\rtimes\{\pm1\}$};
      \node at (3,-3) [name=f,anchor=west] {$Aut_{var}(T_\CC)$};
      \node at (0,-4) [name=g] {$(x,1)$};
      \node at (3,-4) [name=h,anchor=west] {${\exp(x)}$};
      \node at (0,-5) [name=i] {$(x,-1)$};
      \node at (3,-5) [name=j,anchor=west] {${\exp(x)}\circ inv$.};

      \draw[|->] (a) -- (b);
      \draw[->] (c) -- (d);
      \draw[->] (e) -- (f);
      \draw[->] (c) -- node[midway,right] {\small{$\theta'_\bbC$}} (e);
      \draw[->] (3.7,-1.3) -- node [midway,right] {\small{$1$}} (3.7,-2.7);
      \draw[|->] (g) -- (h);
      \draw[|->] (i) -- (j);
    \end{tikzpicture}
  \end{center}
  Here we view the units in $\bbC[\widehat T]$ as a crossed module over $\Aut_{var}(T_\bbC)$ via the trivial map, and we have
  identified $T_\bbC$ with the subgroup of
  $Aut_{var}(T_\bbC)$ given by translations. We use the
  notation $t$ for multiplication with the element $t$ and write $inv$
  for the automorphism $s\mapsto s\inv$.
\end{Def}
The exponential map $\exp$ used in the definition of $r_{bas}$ is an analytic map that does not have an algebraic counterpart.
It is therefore natural to consider the analytic version of the basic representation. For this, we view $T_\bbC$ as a complex analytic variety and replace the algebraic sections $\bbC[\widehat T]=\Gamma\mO_{T_\bbC}$ with the holomorphic sections $\Gamma\mO_{T_\bbC}^{an}$ and $Aut_{var}(T_\bbC)$ with the holomorphic automorphisms $Aut^{an}_{var}(T_\bbC)$. This results in the definition of $r_{bas}^{an}$, the analytic basic representation of $\Theta_\bbC'$.

Given integers $k$ and $n$, we define the representation of level $k$ and weight $n$,
  \begin{center}
  \begin{tikzpicture}
    \node at (-1.5,-2) {$r_{k,n}\negmedspace:$};
    \node at (0,0) [name=a] {$(m,z)$};
    \node at (3,0) [name=b,anchor=west] {$z^{kn}e^{2\pi ikJ^\sharp(m)}$};
    \node at (0,-1) [name=c] {$\Lv\times {\bbC^\times}$};
    \node at (3,-1) [name=d,anchor=west] {$\CC[\,\widehat T\,]^\times$};
    \node at (0,-3) [name=e] {$\fh\rtimes\{\pm1\}$};
    \node at (3,-3) [name=f,anchor=west] {$Aut_{var}(T_\CC)$};
    \node at (0,-4) [name=g] {$(x,1)$};
    \node at (3,-4) [name=h,anchor=west] {${\exp(nx)}$};
    \node at (0,-5) [name=i] {$(x,-1)$};
    \node at (3,-5) [name=j,anchor=west] {${\exp(nx)}\circ inv$.};

    \draw[|->] (a) -- (b);
    \draw[->] (c) -- (d);
    \draw[->] (e) -- (f);
    \draw[->] (c) -- node [midway,right] {$\theta'_\bbC$} (e);
    \draw[->] (3.7,-1.3) -- (3.7,-2.7);
    \draw[|->] (g) -- (h);
    \draw[|->] (i) -- (j);
  \end{tikzpicture}
\end{center}
It, too, has an analytic counterpart $r_{k,n}^{an}$. The basic representation is the special case $r_{bas}=r_{1,1}$.
The strict categorical group $\mA$ associated to the crossed module
$$
  \bbC[\widehat T]\xrightarrow{\,\,\,\,1\,\,\,\,}\Aut_{Var}(T_\bbC)
$$
acts, in a natural manner, on the $\bbC$-linear category $Coh(T_\bbC)$, and similarly in the analytic setting.
Thus, we may interpret $r_{k,n}$ as mapping to a categorical group of $\bbC$-linear functors and linear natural transformations,
  \[
    \begin{tikzpicture}
      \node at (0,0) [name=a]
        {${\varrho_{k,n}}:\,{\mT_\bbC\rtimes\{\pm1\}}$};       
      \node at (4,0) [name=b] {$\mG\mL(\mV)$}; 
      \draw[->] (a) -- (b);
    \end{tikzpicture}
  \]
with $\mV=Coh(T_\bbC)$.
\label{page:rho_kn}  
Restricted to $\mT_\bbC$, the basic representation $\varrho_{bas}=\varrho_{1,1}$
sends
the object $x$ to the functor  
$$
  \longmap{\exp(x)^*}{Coh(T_\CC)}{Coh(T_\CC)},
$$
and it sends the arrow $x\xrightarrow{\,\,\,\,z\,\,\,\,}x+m$ to the
natural transformation
$$
  \varrho_{bas}\(x\xrightarrow{\,\,\,\,z\,\,\,\,}x+m\) \,=\,
  z\exp(x)^*e^{2\pi iJ^\sharp(m)}.
$$  
We do have some freedom in the choice of interpretation of $r_{k,n}$ as action on a category. 
\subsection{Multiplicative bundle gerbe point of view}
\label{sec:multiplicative_gerbe}
\label{sec:connections}
We will write
\[
  T^* = Hom(\Lv,U(1)) \cong \Lambda\tensor U(1) \cong \ft^*/\Lambda
\]
for the Pontryagin dual of $\Lv$ and refer to $T^*$ as the dual torus.
The complexification
\[
  T_\bbC^* = Hom(\Lv, \bbC^\times) = \Lambda\tensor \bbC^\times
\]
may then be viewed as the (analytification of the) complex algebraic
variety
\[
  T_\bbC^* = spec\(\bbC[\check T]\).
\]
We have the Poincar\'e line bundle 
\[
  \mP = \fh\times T_\bbC^*\times \bbC \, / \, \sim
\]
over $T_\bbC\times
T^*_\bbC$, 
with $\sim$ defined by
\[ 
  (x,\eta,z)\sim(x+m,\eta,\eta(m)z).
\]
Writing
\begin{eqnarray*}
  J^\flat_{U(1)} = J^\flat\tensor U(1),
\end{eqnarray*}
we let 
\[
  \mathbf L = \mathbf L^J = \(\id_T\times
  J^{\,\flat}_{U(1)}\)^*\(\mP\at{T\times T^*}\)
\]
be the pull-back of $\mP$ to $T\times T$.
For $t\in T$, we let $\bbL_t$ be the line bundle
  \[
    \bbL_t\, =\,\mathbf L^J_{t\times T}. 
  \]
This bundle comes equipped with
canonical isomorphisms
\[
  \varphi_x\negmedspace : T\times \bbC \,\cong\, \bbL_t,
\]
for any choice of $x$ with $\exp(x)=t$, satisfying
\[
  \varphi_{x+m}\cdot e^{2\pi iJ^\sharp(m)} = \varphi_x.
\]
The bundle $\mathbf L$
was used in \cite{Ganter18}
to equip the trivial bundle gerbe 
over $T$ with a multiplicative structure.
The categorical group corresponding to this multiplicative bundle
gerbe is equivalent to the categorical torus classified by the even bilinear form $I$.
We may think of this as follows. Let $Lines_\bbC$ be the category of $1$-dimensional $\bbC$-vector spaces, and equip 
\[
  \mI = T\times Lines_\bbC,
\]
with the multiplication
\[
  (s,L)\bullet(t,L') = (st,L\tensor L'\tensor \mathbf L_{s,t})
\]
(canonical associators and unit $(1,\bbC)$).\

In this setting the basic representation takes the following shape. 
Let $\mV$ be the category of complex vector bundles 
  on $T$.
  The monoidal functor
  \[
    R_{bas}\negmedspace : (\mI,\bullet) \longrightarrow \(\mG\mL\(\mV\),\circ\)
  \]
  sends the object $(t,L)$ to the endofunctor
  \[
    t_*\circ\(\bbL_{t}\tensor L\tensor -\)
  \] 
  and the pair $((s.L),(t,L'))$ to the natural isomorphism  induced by
  \begin{eqnarray*}
    s_*\(\bbL_{s}\tensor t_* \(\bbL_t \tensor -\)\) & = &
                              s_*t_*\(t^*\bbL_s\tensor\(\bbL_t\tensor
                                                        -\)\) \\[+4pt]
                   & = & (st)_*\(\bbL_{st}\tensor\(\mathbf
                         L^J_{s,t}\tensor -\)\).
  \end{eqnarray*}
This is related to the basic representation $r_{bas}$ of $\mT$ by the diagram
  \begin{center}
    \begin{tikzpicture}
      \node at (0,0) [name=T1] {$\mT$};
      \node at (0,-2) [name=T2] {$\mI$};
      \node at (4,0) [name=A1] {$\mA$};
      \node at (4,-2) [name=A2] {$\mA ut(\mV)$,};
      \draw[->] (T1) -- node [left] {$E$} (T2);
      \draw[->] (T1) -- node [above] {$r_{bas}$} (A1);
      \draw[->] (A1) -- (A2);
      \draw[->] (T2) -- node [below] {$\,\,R_{bas}$} (A2);
      \draw[->, double equal sign distance] (2.7,-.7) -- node [midway,
      above] {$\phi\,\,$} (1.3,-1.3);
    \end{tikzpicture}
\end{center}
 where $E$ takes the object $x$ to $(\exp(x),\bbC)$ and $\phi$ is the
  (monoidal) natural transformation
  \[
    \phi\negmedspace : \varrho_{bas}\,\Longrightarrow\, R_{bas}\circ E
  \]
  given by \[\phi_x\,=\,\exp(x)_*\varphi_x.\]
Finally, the bundle $\boldsymbol L$ inherits from $\mP\at{T\times T^*}$ a connection with holonomy
\[
  Hol_{(\mathbf L,\nabla)} (\bar f,\bar g) = e^{2\pi i
    \(\int_0^1J(\dot f (t),g(t))dt\)-J(f(0),\Delta_g)},
\]
where $\Delta_g=g(1)-g(0)$.
Compare this to \cite{Ganter18}, but with slightly different conventions.
This offers various ways one could hope to reinvent our representations, using complex vector bundles with connection or flat connection in the place of coherent
sheaves.
\subsection{The centre of a categorical torus}
\begin{Lem}
  The centre of the categorical torus $\mT$ is the strict categorical group associated to the crossed module
  \begin{eqnarray*}
    \zeta\negmedspace: \Lv\times U(1) & \longrightarrow & \ft\times\Lambda\\
    (m,z) & \longrightarrow & (m,-I^\sharp(m)),
  \end{eqnarray*}
  with the action given by $$(m,z)^{(x,\lambda)} = \(m,ze^{2\pi i J(m,x)}\).$$
\end{Lem}
\begin{proof}
  Recall that 
  the Pontryagin dual of $\ft$ is identified with the dual vector space via composition with the exponential map,
  \begin{eqnarray*}
    \ft^*&\cong& Hom(\ft,U(1))\\
    \eta & \mapsto & e^{2\pi i\eta(-)},
  \end{eqnarray*}
  see for instance \cite[VII 9.11]{ConwayJB96}.
  It follows that we have
  $$Z_0(\Theta) \ \cong\ \ft\times_{T^*}\ft^*,$$
  where we have made the identification $T^* \cong Hom(\Lv,U(1))$,
  and the map from $\ft$ to $T^*$ sends $x$ to $e^{2\pi i J^b(x)}$.
  The structure map 
  \begin{eqnarray*}
    \Lv\times U(1) & \longrightarrow & \ft\times_{T^*}\ft^*\\
    (m,z) & \longmapsto & (m,-J^\sharp(m)),
  \end{eqnarray*}
  is identified with the map $\zeta$ in the lemma using the isomorphism
  \begin{eqnarray*}
    \ft\times\Lambda & \xrightarrow{\ \cong\ }& \ft\times_{T^*}\ft^*\\
    (x,\lambda) & \longmapsto & (x,J^b(x)+\lambda).
  \end{eqnarray*}
\end{proof}
\begin{Cor}
  In the case where $\mG=\mT$ is the categorical torus, the long exact sequence \eqref{eq:seven-terms} takes the shape
  \[
    \begin{tikzcd}
       1\ar[r]& U(1) \ar[r,"\cong"]&
      U(1)\ar[r]& 1
    \end{tikzcd}
   \]
   \[
   \begin{tikzcd}
    1 \ar[r] & \widehat{T}\ar[r] 
      &\pi_0(\mZ(\Theta))\ar[r]&
      T\ar[r]& 1
    \end{tikzcd}
    \]
    \[
    \begin{tikzcd}
      1\ar[r] & 
      \pi_0\mA ut^+(\mT)\ar[r,"\cong"]&
      \pi_0\mO ut^+(\mT) \ar[r]& 1.
    \end{tikzcd}
  \]
\end{Cor}
  In the self-dual case, i.e., when $I^\sharp$ is an isomorphism, we have a canonical isomorphism $\ft\cong\pi_0(\mZ(\Theta))$, and the short exact sequence in the
  middle row is identified with the usual short exact sequence $$0 \longrightarrow\Lv\longrightarrow  \ft\xrightarrow{\ \exp\ } T\longrightarrow 1.$$
\section{Centralisers and normalisers}\label{sec:centralisers}
Let $k$ be a field, and let
$$
  \varrho\negmedspace : G\longrightarrow GL(V)
$$
be a $k$-linear representation of a finite group. Then we have the centraliser and the normaliser of $\varrho$, defined as
$$
  C(\varrho) = \{h\in GL(V)\mid \varrho = c_h\circ \varrho\}
$$
and
$$
  N(\varrho) =   N(\varrho) = \{(f,h)\in Aut(G)\times GL(V)\mid \varrho\circ f=h\varrho(-) h\inv\}.
$$
These sit inside a commuting diagram with exact rows
\[
  \begin{tikzcd}
    1\ar[r]& Z(G)\ar[d]  \ar[r] &G\ar[d]\ar[r] & Inn(G)\ar[r] \ar[d]& 1\\
    1\ar[r]& C(\varrho)  \ar[r] &N(\varrho) \ar[r] & Aut(G),
  \end{tikzcd}
\]
where $Z(G)$ denotes the centre of $G$.
If $\varrho$ is irreducible, then Schur's lemma implies
$$C(\varrho)=k^\times.$$
Let $(\Lambda,I)$ be the Leech lattice equipped with its even symmetric bilinear form, and consider the extraspecial $2$-group $\widetilde{\Lambda/2\Lambda}$, with commutator given by $$e_\alpha e_\beta = (-1)^{I(\alpha,\beta)} e_\beta e_\alpha,\qquad\alpha,\beta\in\Lambda/2\Lambda.$$
Let $\phi$ be the quadratic form
\begin{eqnarray*}  
  \phi\negmedspace :\Lambda\tensor_\bbZ\bbF_2&\longrightarrow &\bbF_2\\
  \alpha & \longmapsto & \frac12 I(\alpha,\alpha),
\end{eqnarray*}
and $O(\phi)$ its isometry group.
The construction of the monster
begins by considering the action of $\widetilde{\Lambda/2\Lambda}$
on its unique centre-faithful irreducible representation $\varrho$. This representation is a $2^{12}$-dimensional real representation, whose centraliser and normaliser fit into the diagram
\[
  \begin{tikzcd}
   &&1\ar[d]&1\ar[d]&\\
     1\ar[r]& \{\pm1\}\ar[d,equals]  \ar[r]
     &\widetilde{\Lambda/2\Lambda} \ar[d]\ar[r] &
     \Lambda/2\Lambda\ar[r] \ar[d]& 1\\
    1\ar[r]& C(\varrho)  \ar[r] &N(\varrho) \ar[r]\ar[d] & Aut(\widetilde{\Lambda/2\Lambda})\ar[r]\ar[d] & 1\\
    & &O(\phi)\ar[d]\ar[r,equals]& Out(\widetilde{\Lambda/2\Lambda})\ar[d] & \\
          & &1& 1& \\
  \end{tikzcd}
\]
\cite[10.4.15]{FrenkelLepowskyMeurman88}.
The exactness on the right in the second row follows from the uniqueness of $\varrho$.
The group $O(\phi)$ receives a map from $Co_1$ and hence $Co_0$, pulling back the middle column to an \label{page:Ctilde} extension 
$$
  \begin{tikzcd}
    1\ar[r]& \widetilde{\Lambda/2\Lambda}  \ar[r] &\widetilde C\ar[r]
    & Co_0\ar[r]& 1.
  \end{tikzcd}
$$
One then identifies the central elements $-1\in\widetilde{\Lambda/2\Lambda}$ and $-\id\in Co_0$ to obtain the group
$$
  C = \widetilde C\,/\,(-1)\sim(-\id).
$$
This is an extension of $Co_1$ by the extraspecial 2-group  $\widetilde{\Lambda/2\Lambda}$.  
The Monster can be generated by $C$ and one more element, the triality symmetry.
\subsection{The centraliser of a crossed module homomorphism}

Given a strict homomorphism of crossed modules
  \[
    \begin{tikzcd}
      G_1\ar[r,"f_1"] \ar[d,swap,"\psi"] & H_1\ar[d,"\phi"]\\
      G_0\ar[r,"f_0"] & H_0,
    \end{tikzcd}
  \]
  as in Definition \ref{def:homomorphism}, we define the group $C_0(f)$ as follows. Elements of $C_0(f)$ are pairs $(h,\chi)$, with $h\in H_0$ and $\chi\negmedspace:G_0\longrightarrow H_1$, such that
\begin{eqnarray*}
  \ul{\chi(s)} & \ =\ & f_0(s\inv)h\inv f_0(s) h,\\
  \chi(\ul\alpha) & = & f_1(\alpha\inv) f_1(\alpha)^h,\\
  \chi(st) &=& \chi(s)^{f_0(t)}\chi(t).
\end{eqnarray*}
In the Lie setting, we require the map $\chi$ to be analytic.
The multiplication in $C_0(f)$ is given by
$$(h,\chi)(k,\sigma)\ =\ (hk, s\mapsto \sigma(s)\chi(s)^k).$$
\begin{Def}
  In the above situation, the centraliser of $f$ is the crossed module $C(f)$ given by
  \begin{eqnarray*}
    \gamma\negmedspace : H_1 & \longrightarrow & C_0(f)\\
    \zeta & \longmapsto & \(\ul\zeta,s\mapsto \(\zeta^{f_0(s)}\)\inv\zeta\)
  \end{eqnarray*}
  where $C_0(f)$ acts on $H_1$ via
  $\zeta^{(h,\chi)} = \zeta^h$.
\end{Def}
A computation closely following \cite{pirashvili2023centre} shows that $C(f)$ is indeed a crossed module. The centre of a crossed module is a special case of this construction, $$Z(\Psi) = C\(\id_\Psi\),$$ and
there is a crossed module homomorphism
\begin{eqnarray*}
  Z(\Psi) & \longrightarrow & C(f)\\
  (s,\xi) & \longmapsto & (f_0(x),f_1\circ\xi)\\
  \alpha & \longmapsto & f_1(\alpha).
\end{eqnarray*}
\subsection{The centraliser of the basic representation of a categorical torus}
We are interested in the centraliser $C(r_{bas}^{an})$, where $r_{bas}^{an}$ is the crossed module homomorphism of Definition \ref{def:basic_rep}, encoding the analytic basic representation of $\mT_\bbC\rtimes\{\pm1\}$.
\begin{Def}
  We will write $E$ for the crossed module
  \begin{eqnarray*}
    \eta\negmedspace : \Lv\times\bbC^\times & \longrightarrow & \widetilde\Lv\\
    (m,z) & \longmapsto & 2m,
  \end{eqnarray*}
  where the central extension is defined by the $2$-cocycle $J^t$ modulo $2$, 
  \begin{eqnarray*}
    \Lv\times\Lv & \longrightarrow & \bbZ/2\bbZ\\
    (m,n) & \longmapsto &  [J(n,m)], 
  \end{eqnarray*}
  and acts on $\Lv\times\bbC^\times$ via $(m,z)^{(n,\iota)} = (m,(-1)^{J(n,m)}z)$.
\end{Def}
The strict categorical group associated to $E$ is the extraspecial categorical $2$-group $$\mS(E)=\mT_\bbC^{\{\pm1\}}$$  
with
\begin{eqnarray*}
  \pi_0(E)\ \ \cong\ \ \widetilde{T[2]}& \cong& \widetilde{\Lv/2\Lv}
\end{eqnarray*}
and $\pi_1(E)=\bbC^\times$, which 
turned up in \cite[Theorem 6.3]{Ganter18} as the categorical fixed points of the action of $\{\pm1\}$ on the complexification of $\mT\simeq\mT_{J^t}$. This is the same action as the one in \eqref{eq:iota} above.
We find it convenient to use additive notation for the centre of $\widetilde{\Lv}$.
\begin{Thm}
  We have an equivalence of crossed modules $F\negmedspace :E\longrightarrow C_0\(r_{bas}^{an}\)$.  
\end{Thm}
\begin{proof}
Define $F$ to be the strict crossed-module homomorphism
  \[
    \begin{tikzcd}
      (m,z) \ar[rr,mapsto]\ar[dd,mapsto] &[-.75cm] &[+1cm] ze^{-2\pi i J^b(m)}&[-1cm]\\[-.5cm]
      & \Lv\times\bbC^\times\ar[r,"F_1"]\ar[d,"\eta"]
      &\(\Gamma\mO_{T_\bbC}^{an}\)^\times\ar[d,"\gamma"] & \zeta\ar[d,mapsto]\\[+1cm]
      (2m,0) & \widetilde\Lv\ar[r,"F_0"]& C_0\(r_{bas}^{an}\) & \(\id,\chi_\zeta\)\\[-.5cm]
      &(n,\iota) \ar[r,mapsto]& \(\exp\(\frac n2\),\chi_{(n,\iota)}\) 
    \end{tikzcd}
  \]
  where
  \[
    \chi_{(n,\iota)}(x,\epsilon)  = \epsilon^\iota e^{\pi i\(J(x,n)+(\epsilon-1)J^b(n)\)}
  \]
  and
  \[
    \chi_\zeta(x,\epsilon) (t) = \frac{\zeta(t)}{\zeta(\exp(x)\cdot t^\epsilon)}.
  \]
  If $\zeta=ze^{2\pi i\lambda}$ then $$\chi_\zeta(x,\epsilon) = e^{-2\pi i(\lambda(x)+(\epsilon-1)\lambda)},$$
  so the diagram commutes.
  The section $\zeta\negmedspace:T_\bbC\longrightarrow \bbC^\times$ is in  $\pi_1\(C\(r_{bas}^{an}\)\)$ if and only if $\zeta$ satisfies
  $$\zeta(\exp(x)t^\epsilon)=\zeta(t)$$
  for all $(x,\epsilon)\in\fh\rtimes\{\pm1\}$ and
  $t\in T_\bbC$. Taking $t=1$, we see that if $\zeta$ is an element of $\pi_1\(C\(r_{bas}^{an}\)\)$, then $\zeta$ is constant on $T_\bbC$. Hence $\pi_1(F)$ is an isomorphism.
  To see that $\pi_0(F)$ is injective, assume we are given $(n,\iota)\in\widetilde\Lv$ and $\zeta\in \(\Gamma\mO_{T_\bbC}^{an}\)^\times$ such that  
  $\exp\(\frac n2\)=1$ and $\chi_{(n,\iota)} = \chi_\zeta$.
  Then $n\in2\Lv$, and we have
  \[
     {\zeta(\exp(x)\cdot t^\epsilon)} = \zeta(t)\cdot \epsilon^\iota e^{-\pi i\(J(x,n)+(\epsilon-1)J^b(n)\)}(t)
  \]
  for all $(x,\epsilon)\in\fh\rtimes\{\pm1\}$ and $t\in T_\bbC$. 
  Setting $t=1$ and $z=\zeta(1)$, we otain
  \[
     {\zeta(\exp(x))} = z \cdot \epsilon^\iota e^{-2\pi iJ\(x,\frac{n}2\)}.
  \]
  The left-hand side is independent of $\epsilon$, and we conclude that $\iota=0$ and $\zeta=z e^{-2\pi i J^b\(\frac{n}2\)}$.
  It remains to show surjectivity of $\pi_0(F)$.
  Assume that we are given an automorphism $f$ of $T_\bbC$ commuting with multiplication by $\exp(x)$ for all $x\in\fh$. Then 
  $$f(\exp(x)) = \exp(x)\cdot f(1),$$ so $f$ is given by multiplication with $f(1)$. If, in addition, we assume that $f(t)\inv=f(t\inv)$ holds for all $t\in T_\bbC$, it
  follows that $f(1)\in T[2]$ is an element of order $2$. Choose $n\in\Lv$ such that $f(1)=\exp\(\frac n2\)$, and let
  $$\chi\negmedspace:\fh\rtimes\{\pm1\}\longrightarrow \(\Gamma\mO_{T_\bbC}^{an}\)^\times$$
  be such that $\(\exp\(\frac n2\),\chi\)$ is an element of $C_0\(r_{bas}^{an}\)$. If $m\in\Lv$, then $\chi$ sends the element $(m,1)$ to the constant function
  $(-1)^{J(m,n)}$. Hence 
  \begin{eqnarray*}
    \zeta\negmedspace: T_\bbC & \longrightarrow & \bbC^\times\\
    t &\longmapsto& \frac{e^{\pi iJ(y,n)}}{\chi(y,1)(1)},
  \end{eqnarray*}
  where $t=\exp(y)$, is a well-defined element of $\(\Gamma\mO_{T_\bbC}^{an}\)^\times$. Further,
  we have $$\chi(0,-1)\(t\inv\)\cdot \chi(0,-1)(t) = \chi(0,1)(t) = 1,$$
  and hence
  $\chi(0,\epsilon)(1) = \epsilon^\iota$ for some $\iota\in\bbZ/2\bbZ$, and $$\chi(x,\epsilon)(1) = \chi(x,1)\(1\)\cdot\epsilon^\iota.$$
  Letting $t=\exp(y)$, we compute
  \begin{eqnarray*}
    \frac{\zeta(t)}{\zeta(\exp(x)\cdot t^\epsilon)}
    & = & \frac{\chi(x+\epsilon y,1)}{\chi(y,1)}(1)\cdot e^{\pi iJ(-x-\epsilon y+y,n)}\\
                                                   & = & \frac{\chi(x+\epsilon y,\epsilon)}{\chi(y,1)}(1)\cdot \epsilon^\iota e^{-\pi iJ(x+\epsilon y-y,n)}\\
                                                   & = & \frac{\chi(x,\epsilon)}{\chi_{(n,\iota)}}(t).
  \end{eqnarray*}
  It follows that the classes of $\(\exp\(\frac n2\),\chi\)$ and $\(\exp\(\frac n2\),\chi_{n,\iota}\)$ in $\pi_0\(C\(r_{bas}^{an}\)\)$ agree, and we have shown that $\pi_0(F)$ is surjective.
\end{proof}

\section{Character theory}

The character theory of categorical representations was developed in 
\cite{Bartlett11} and
\cite{GanterKapranov08} and has recently seen a revival in a slightly diffferent setting: the theory relies on a notion of categorical trace,
and there are two such notions in use. These are different, but related; 
see \cite{Willerton2Traces} for a comparison and 
\cite{gaitsgory2020local} for a recent application of categorical traces.
We will be working with the definition
\[
  \ttr (F) = 2Hom(\id,F)
\]
for the categorical trace of a 1-endomorphism $F$.
By the {inertia groupoid} 
\begin{eqnarray*}
   G/\negmedspace/ G^{conj}  & \cong & Fun(\bbB \bbZ, \bbB G)
\end{eqnarray*}
we will mean the action groupoid of $G$ acting on itself by conjugation. Characters of classical representations are invariant functions on
$G/\negmedspace/ G^{conj} $.
If a finite group $G$ acts either by
functors on a $k$-linear category, or more generally, by 1-morphisms 
in a $k$-linear bicategory, we speak of a $k$-linear categorical representation of $G$.
In this
context, there are two levels of character theory.
Applying the categorical trace, one obtains
the categorical character $X_\varrho$. This is a $k$-vector bundle over $G/\negmedspace/ G^{conj}$, meaning that one has
$k$-vector spaces $$X_\varrho(g) = \ttr(\varrho(g)),$$ and compatible isomorphisms
$$\psi_{g,h}\negmedspace: X_\varrho(g)\cong X_\varrho(h g h\inv).$$
If $g$ and $h$ commute, then $\psi_{g,h}$ is an endomorphism of $\ttr(g)$, and in many interesting settings one can take the trace again.
This results in the definition of 2-character $$\chi(g,h)=tr(\psi_{g,h}),$$ which is a function on pairs of commuting elements of $G$, invariant under simultaneous conjugation. In other words, the 2-character is
 an invariant function on the interated inertia groupoid $Fun(\bbB\bbZ^2, \bbB G)$.

Consider next an action of a finite categorical
group $\mG$. This may be viewed as a projective 
categorical representation of the group
$G=\pi_0(\mG)$. 
By a categorical class function on $\mG$ we shall mean a bundle on the {inertia
  groupoid} of $\mG$,
\begin{eqnarray*}
  \mG/\negmedspace/\mG^{conj} & := & Bifun(\bbB\bbZ,\bbB \mG)\,\slash\, \text{2-isos}.
\end{eqnarray*}
It was shown in {\cite{GanterUsher16}}\label{thm:Usher} that
the categorical character of a $k$-linear representation
$\varrho$ of a finite categorical group $\mG$ is a
categorical class function on $\mG$.
If $\mG$ is classified by $[\alpha]\in H^3(G;k^\times)$, then $\mG/\negmedspace/\mG^{conj}$ is equivalent to the central extension of
$$G/\negmedspace/G^{conj}\simeq\coprod_{[g]}1/\negmedspace/\widetilde C_g$$ classified by the 
transgression of $\alpha$. 
This consists of central extensions $\widetilde C_g$ of the centralisers by $k^\times$, one for each conjugacy class $[g]$ in $G$.
The categorical character of $\varrho$ is thus
a module over the $\alpha$-twisted Drinfeld double of $G$ in the sense of \cite{Willerton08}.
When $G=M$ is the Monster and $\alpha=\amoon$ is the Moonshine anomaly, this resembles the formalism governing generalised Moonshine.

  Turning our attention to the torus, the inertia groupoid of the commutative group $T$ is easily described: it has objects $T$ and arrows $T\times T$ with source and target both given by the projection to the first factor, the fibre $\{t\}\times T$ being the centraliser of $t$ in $T$. There are no other arrows, since everything commutes. In light of the picture for finite categorical groups, we expect $\mT/\negmedspace/\mT^{conj}$ to be a central extension of $T/\negmedspace/T^{conj}$ by $U(1)$.

  This means that we expect to find each centraliser $T=C_t$ replaced by a central extension by $U(1)$. As individual groups, these are trivial central extensions, so we have non-canonical isomorphisms $\widetilde C_t\cong T\times U(1)$.
  We will, however, find that with varying $t$ the global topology of the inertia groupoid of $\mT$ is non-trivial. 
\subsection{Inertia 2-groupoids}
For a brief introduction to the language of bicategories we refer the reader to
\cite{leinster1998basic}, noting that we will follow the opposite convention for pseudo-natural transformations,
found, for instance, in \cite{GordonPowerStreet95}.
  Let $\mG$ be a categorical group.
  By the inertia 2-groupoid of
  $\mG$ we will mean the bifunctor 2-groupoid
 $
    \Bicat(\bbB\bbZ,\bbB\mG).
  $
Its objects may be thought of as strong monoidal functors from $\bbZ$ to
$\mG$. The small inertia 2-groupoid $\Bicat^{strict}(\bbB\bbZ,\bbB\mG)$ is defined as the full sub-2-groupoid whose objects are strict monoidal functors. 
We will see that the inclusion of the small inertia 2-groupoid inside the full inertia 2-groupoid for $\mT$ is an equivalence of Lie 2-groupoids.
The inertia groupoid $\mG/\negmedspace/\mG^{conj}$ is obtained from the inertia 2-groupoid by taking 2-isomorphism classes of $1$-arrows.
The goal of this section is to describe $\mT/\negmedspace/\mT^{conj}$.
%
\begin{Lem}
\label{lem:HZ}
  Let $M$ be a trivial $\bbZ$-module, and let $N$ be a submodule of
  $M$.
  \begin{enumerate}[(i)]
  \item Any 1-cochain $\longmap\mu\bbZ M$ is uniquely determined by
    $d\mu$ and $\mu(1)$.
  \item The 1-cochain $\mu$ as above takes values in $N$ if and only
    if $d\mu$ does and $\mu(1)\in N$.
  \item For any 2-cocycle $\longmap\gamma{\bbZ\times\bbZ}M $ and integer
    $a\in\bbZ$, we have 
    $$
      \gamma(0,a) \,\,=\,\,\gamma(0,0)\,\,=\,\,\gamma(a,0).
    $$   
  \end{enumerate}
\end{Lem}
\begin{proof}
  (i) If $\mu_1$ and $\mu_2$ are 1-cochains with $d\mu_1=d\mu_2$, then
  their difference is a group homomorphism $\longmap{\mu_1-\mu_2}\bbZ
  M$. If in addition $\mu_1(1)=\mu_2(1)$, it follows that
  $\mu_1-\mu_2=0$.

  \smallskip
  \noindent
  (ii) Assume that $d\mu$ takes values in $N$. Since $H^2(\bbZ,N) = 0$,
  there exists a 1-cochain $\longmap\nu\bbZ N$ with $d\nu=d\mu$. Using
  (i), we deduce
  \begin{eqnarray*}
    \mu(a)  & = &  \nu(a) + a(\mu(1)-\nu(1)).
  \end{eqnarray*}
  So, if $\mu(1)\in N$, it follows that $\mu$ takes values in $N$.

  \smallskip
  \noindent
  (iii) As in (ii), choose any 1-cochain $\longmap\alpha\bbZ M$ with
  $d\alpha=\gamma$. Let $a$ be an integer. Then
  $$
    \gamma(0,a) \,\,=\,\,\alpha(0)\,\,=\,\,\gamma(a,0).
  $$   
\end{proof}
\begin{Thm}\label{thm:full_inertia_2-groupoid}
  The inertia 2-groupoid $\Bicat(\bbB\bbZ,\bbB\mT)$ is isomorphic to the Lie
  2-groupoid with 
  \begin{description}
  \item[objects] triples $(x,\gamma,c)$ with $x\in \ft$, and
    \begin{eqnarray*}
      \gamma\negmedspace : \bbZ\times\bbZ & \longrightarrow & \Lv \\
      c\negmedspace : \bbZ\times\bbZ & \longrightarrow & U(1) 
    \end{eqnarray*}
    2-cocycles;
  \item[1-arrows] six-tuples $(x,\gamma,c,y,\mu,w)$, with source and
    target as follows
    \begin{center}
      \begin{tikzpicture}
        \node at (0,0) [name=a] {$(x,\gamma,c)$};
        \node at (5,0) [name=b] {$\(x+\mu(1),\gamma+d\mu, c'\)$.};
        \draw [->] (a) -- (b) node [midway, above] {$(y,\mu,w)$};
      \end{tikzpicture}  
    \end{center}
    Here we have $(x,\gamma,c)$ as above, $y\in\ft$,
    \begin{eqnarray*}
      \mu\negmedspace : \bbZ & \longrightarrow & \Lv \\
      w\negmedspace : \bbZ & \longrightarrow & U(1) 
    \end{eqnarray*}
    1-cochains, and
    \begin{eqnarray*}
      c'(a,b) & = & c(a,b)\cdot dw(a,b)\cdot \(e^{2\pi i J\(\mu(1),x\)}\)^{ab};
    \end{eqnarray*}
  \item[2-arrows] eight-tuples $(x,\gamma,c,y,\mu,w,n,u)$, viewd as 2-arrow
    \begin{center}
      \begin{tikzpicture}
        \node at (2.3,0) {$\Downarrow u$};
        \node at (0,0) [name=a] {$(x,\gamma,c)$};
        \node at (4,0) [name=b, anchor = west] {$\(x+\mu(1),\gamma+d\mu, c'\)$};
        \draw [->,bend left=20] (a) to  node
        [midway, above] {$(y,\mu,w)$} (b.north west); 
        \draw [->,bend right=20] (a) to node
        [midway, below] {$(y+n,\mu,w')$} (b.south west); 
      \end{tikzpicture}  
    \end{center}
    with $n\in\Lv$, $u\in U(1)$, and
    \begin{eqnarray*}
      w'(a) & = & w(a)\cdot \(e^{2\pi i I(x,n)}\)^a;
    \end{eqnarray*}
  \item[Horizontal composition] the following diagram commutes strictly
    \begin{center}
      \begin{tikzpicture}
        \node at (-1,0) [name=a] {$(x,\gamma,c)$};
        \node at (3,1.5) [name=b] {$\(x+\mu(1),\gamma+d\mu, c'\)$};
        \node at (7.5,0) [name=c] {$\(x+(\mu+\mu')(1),\gamma+d(\mu+\mu'), c''\)$};
        \draw [->] (a) -- (b) node [midway, anchor=south east] {$(y,\mu,w)$};
        \draw [->] (b) -- (c) node [midway, anchor=south west]
        {$(y',\mu',w')$};
        \draw [->] (a) -- (c) node [midway, below]
        {$(y+y',\mu+\mu',w\cdot w')$};
      \end{tikzpicture}  
    \end{center}
    where 
    \begin{eqnarray*}
      c''(a,b) & = & c(a,b) \cdot dw''(a,b) \cdot 
          \(e^{2\pi iJ\((\mu+\mu')(1),x\)}\)^{ab}.
    \end{eqnarray*}
  \end{description}
\end{Thm}
\begin{proof}
  {\bf The objects:}
  a strong monoidal functor $\bbZ\longrightarrow \mT$ consists of 
  a sequence $a\longmapsto x_a$, indexed by $\bbZ$, together with
  arrows
    \begin{eqnarray*}
      \varphi_{a,b}\negmedspace : x_a+x_b &
                                            \stackrel\cong\longrightarrow & x_{a+b}
    \end{eqnarray*}
    and
    \begin{eqnarray*}
      \phi_{0}\negmedspace : 0 &
                                 \stackrel\cong\longrightarrow & x_{0},
    \end{eqnarray*}
    satisfying
    \begin{eqnarray}
      \label{eq:hexagon}
      \phi_{a+b,c}\circ (\phi_{a,b}\bullet\id_{x_c}) & = & 
                                                           \phi_{a,b+c}\circ (\id_{x_a}\bullet\,\phi_{b,c}) 
    \end{eqnarray}
    and 
    \begin{equation}
      \label{eq:unit}
      \phi_{0,a}\circ(\phi_0\bullet\id_{x_a})\,\,=\,\,\id_{x_a}\,\,=\,\,
      \phi_{a,0}\circ(\id_{x_a}\bullet\phi_0).
    \end{equation}
  Given such data, we set $x:= x_1$ and $\lambda(a):= x_a-ax$ and $\gamma=d\lambda$. This
  allows us to write $\phi_{a,b}$ and $\phi_0$ in the following form
  %
  %
  \begin{center} 
    \begin{tikzpicture}
      \node at (0,0) [name=a] {$\phi_{a,b}: \quad x_a+x_b$};
      \node at (7,0) [name=b] {$x_a+x_b+\gamma(a,b)$.};
      \draw [->] (a) -- (b) node [midway, above] {$c(a,b)\cdot e^{J(\lambda(a),bx)}$};
    \end{tikzpicture}  
  \end{center}
  and
  \begin{center}
    \begin{tikzpicture}
      \node at (0,0) [name=a] {$\phi_{0}: \quad 0$};
      \node at (4,0) [name=b] {$\lambda(0)$.};
      \draw [->] (a) -- (b) node [midway, above] {$c(0,0)\inv$};
    \end{tikzpicture}  
  \end{center}
  The condition \eqref{eq:hexagon} is satisfied if and only if
  \begin{eqnarray*}
    c\negmedspace : \bbZ\times \bbZ & \longrightarrow & U(1)
  \end{eqnarray*}
  is a 2-cocycle. By Lemma \ref{lem:HZ}, the 1-cochain $\lambda$
  takes values in $\Lv$ and is uniquely determined by $\gamma$ and by
  $\lambda(1)=0$.
  To summarize, the data of a strong monoidal functor
  $\bbZ\longrightarrow \mT$ are determined by the triple $(x,\gamma,c)$,
  where $x\in \ft$, and $\longmap\gamma{\bbZ\times\bbZ}\Lv$ and
  $\longmap c{\bbZ\times\bbZ}U(1)$ are 2-cocycles.
  
  {\bf The 1-arrows:}
  a pseudo-natural transformation from $((x_a)_a,\phi_{a,b},\phi_0)$ to 
  $((x_a')_a,\phi_{a,b}',\phi_0')$, viewed as bifunctors $\bbB\bbZ\longrightarrow\bbB G$, consists of
  an element $y\in \ft$, together with a sequence of arrows
  \begin{eqnarray*}
    \beta_a\negmedspace : y+x_a & \stackrel\cong\longrightarrow & x_a'+y
  \end{eqnarray*}
  satisfying
  \begin{eqnarray}
    \label{eq:octagon}
    \beta_{a+b}\circ \(\id_y\,\bullet \,\,\phi_{a,b}\) & = &
                                                             \(\phi_{a,b}'\,\bullet\,\,\id_y\) \circ \(\id_{x_a'}\,\bullet
                                                             \,\,\beta_{b}\) \circ \(\beta_a\,\bullet\,\,\id_{x_b}\) 
  \end{eqnarray}
  and 
  \begin{eqnarray}
    \label{eq:unit2}
    \beta_0\circ (\id_y\,\bullet \,\,\phi_0) & = & \phi_0'\, \bullet\, \id_y.
  \end{eqnarray}
  Writing $\beta_a$ in the form
  \begin{center}
    \begin{tikzpicture}
      \node at (0,0) [name=a] {$\beta_{a}: \quad y+x_a$};
      \node at (7,0) [name=b] {$x_a+\mu(a)+y$,};
      \draw [->] (a) -- (b) node [midway, above] {$w(a)\cdot e^{-J(x_a+\mu(a),y)}$};
    \end{tikzpicture}  
  \end{center}
  condition \eqref{eq:octagon} spells out to 
  \begin{eqnarray*}
    d\mu & = & d\lambda'-d\lambda, \quad\text{and}\\
    dw(a,b) & = & \frac{c'(a,b)}{c(a,b)}\cdot\(e^{-J(\mu(1),x)}\)^{ab}.
  \end{eqnarray*}
  Let $m=\mu(1)$ and $z=w(1)$. Then $\mu$ and $w$ are uniquely
  determined by their boundaries together with $m$ and $z$.
  To summarize, a 1-arrow in $\Bicat(\bbB\bbZ,\bbB\mT)$ between the monoidal
  functors determined by $(x,\gamma,c)$ and $(x+m,\gamma',c')$ is
  determined by its source and target together with a pair 
  $(y,z) \in\ft\times U(1)$. Any choice of $(y,z)$ is valid.

  {\bf The 2-arrows}
  These are the modifications.
\end{proof}
\begin{Lem}
  The inclusion of the strict inertia 2-groupoid inside the full inertia 2-groupoid of $\mT$ is an equivalence.
\end{Lem}
  \begin{proof}
    The set of 1-arrows from $(x,0,1)$ to
    $(x,\gamma,c)$ is non-empty. Given a choice of $(y,z)$, the arrows 
    $\beta_a$ look as follows:
    \begin{center}
      \begin{tikzpicture}
        \node at (0,0) [name=a] {$\beta_{a}: \quad y+ax$};
        \node at (8,0) [name=b] {$ax+\lambda(a)+y$.};
        \draw [->] (a) -- (b) node [midway, above]
        {$w(a)\cdot e^{-J(ax+\lambda(a),y)}$}; 
      \end{tikzpicture}  
    \end{center}
    Here $\longmap\lambda\bbZ\Lv$ is the 1-cochain with $d\lambda=\gamma$ and
    $\lambda(1)=0$, and $\longmap w\bbZ U(1)$ is the 1-cochain with $dw=c$ and
    $w(1)=z$, as in Lemma \ref{lem:HZ} (i).
  \end{proof}
%
\begin{Cor}\label{cor:small_inertia_2-groupoid}
  The small inertia 2-groupoid of $\mT$ is isomorphic to the Lie
  2-groupoid with 
  \begin{description}
  \item[objects] elements $x\in \ft$;
  \item[1-arrows] quadruples 
    \begin{eqnarray*}
      (x,y,m,w) & \in & \ft\times \ft\times \Lv\times U(1)
    \end{eqnarray*}
    written
    \begin{center}
      \begin{tikzpicture}
        \node at (0,0) [name=a] {$x$};
        \node at (2.5,0) [name=b] {$x+m$.};
        \draw [->] (a) -- (b) node [midway, above] {$(y,w)$};
      \end{tikzpicture}  
    \end{center}
  \item[2-arrows] as follows
    \begin{center}
      \begin{tikzpicture}[scale=1]
        \node at (2,0) {$\Downarrow u$};
        \node at (0,0) [name=a] {$x$};
        \node at (3.8,0) [name=b, anchor = west] {$x+m$};
        \draw [->, bend left=30] (a.north east) to node
        [midway, above] {$(y,w)$}  (b.north west); 
        \draw [->, bend right=30] (a.south east) to node
        [midway, below] {$\(y+n,w\cdot e^{I(x,n)}\)$} (b.south west); 
      \end{tikzpicture}  
    \end{center}
    with $n\in\Lv$ and $u\in U(1)$.
  \item[Horizontal composition] the following diagram commutes strictly
    \begin{center}
      \begin{tikzpicture}
        \node at (0,0) [name=a] {$x$};
        \node at (3.5,1.5) [name=b] {$x+m,$};
        \node at (7,0) [name=c] {$x+m+m'$.};
        \draw [->] (a) -- (b) node [midway, anchor=south east] {$(y,w)$};
        \draw [->] (b) -- (c) node [midway, anchor=south west]
        {$(y',w')$};
        \draw [->] (a) -- (c) node [midway, below]
        {$(y+y',m+m',w\cdot w')$};
      \end{tikzpicture}  
    \end{center}
  \end{description}
\end{Cor}
\begin{proof}
  The full inertia 2-groupoid is described in Theorem \ref{thm:full_inertia_2-groupoid}.
  In the notation of Theorem \ref{thm:full_inertia_2-groupoid}, set
  $$m:=\mu(1)\quad\text{ and }\quad w:=w(1).$$ Then the 1-cocycle $\mu$ and the 1-cochain
  $w$ are uniquely determined by $m$ and $w$.
\end{proof}

\begin{Thm}\label{thm:inertia_groupoid}
  The inertia groupoid $\mT/\negmedspace/\mT^{conj}$ is equivalent 
  to 
  \begin{enumerate}[(i)]
  \item the Lie groupoid with objects $T$ and arrows
  $$
    T\times\ft\times U(1) /\sim
  $$
  with 
  \begin{eqnarray*}
    (t,y+n,w) & \sim & \(t,y,w\cdot e^{-I(x,n)}\),\quad\quad n\in\Lv,
  \end{eqnarray*}
  where $x\in \ft$ is any element with $\exp(x)=t$, source and target
  of $(t,y,w)$ equal $t$, and composition of arrows is
  \begin{eqnarray*}
    [t,y,w]\circ[t,y',w'] & = & [t,y+y',w\cdot w'];
  \end{eqnarray*}
  \item the Lie groupoid with objects $T$ and arrows
  $$
    \mathfrak t \times T \times U(1) /\sim
  $$
  with 
  \begin{eqnarray*}
    (x+m,s,z) & \sim & \(x,s,z\cdot e^{I(m,y))}\),\quad\quad m\in\Lv,
  \end{eqnarray*}
  where $y\in \ft$ is any element with $\exp(y)=s$, source and target
  of $(x,s,z)$ equal $\exp(x)$, and composition of arrows is
  \begin{eqnarray*}
    [x,s,z]\circ[x,s',z'] & = & [x,s\cdot s',z\cdot z'].
  \end{eqnarray*}
  \end{enumerate}
\end{Thm}
\begin{proof}
  Part (i) is immediate from Corollary
  \ref{cor:small_inertia_2-groupoid}. An isomorphism between the
  groupoids described in the two parts is given by
  \begin{eqnarray*}
    \ft\times\ft\times U(1) & \longrightarrow & \ft\times\ft\times
    U(1) \\
    (x,y,w) & \longmapsto & (x,y,z) 
  \end{eqnarray*}
  with $z=w\cdot e^{-I(x,y)}$. This isomorphism is
  $\Lv\times\Lv$ equivariant with respect to the actions described in
  the Theorem.
\end{proof}
If we view $\mT/\negmedspace/\mT^{conj}$ as family of $U(1)$-central extensions of centralisers,
$\widetilde C_t$, parametrized by $t\in T$, then Theorem
\ref{thm:inertia_groupoid}(ii) shows that each $\widetilde C_t$ is
isomorphic to the trivial central extension
\begin{eqnarray*}
  \widetilde C_t & \cong & T\times U(1),
\end{eqnarray*}
but there is no global trivialization over $T$. Instead, the
trivialization for $\widetilde C_t$ depends on a choice of $x\in\ft$
with $\exp(x)=t$. Different choices yield trivializations that differ
by an automorphisms of $T\times U(1)$. More precisely, 
we have a commuting diagram
\begin{center}
  \begin{tikzpicture}[scale=1.5]
    \node at (-.5,0) [name=a] {$\widetilde C_t$};
    \node at (2,.75) [name=b] {$T\times U(1)$};
    \node at (2,-.75) [name=c] {$T\times U(1)$,};

    \draw[->] (a)--(b) node[midway,anchor=south,yshift=.2cm]{$triv_x$};
    \draw[->] (a)--(c) node[near end,anchor=north
    east]{$triv_{x+m}$};
    \draw[->] (b)--(c) node[midway,anchor= west]{};
  \end{tikzpicture}
\end{center}
where the vertical arrow is
obtained by exponentiating the map
\begin{eqnarray*}
  \(
  \begin{matrix}
    \mathds 1&0\\-I^\sharp(m)&1
  \end{matrix}
\): \ft\times\bbR & \longrightarrow & \ft\times\bbR.
\end{eqnarray*}
\subsection{Class functions on categorical tori}\label{sec:class_functions}
Consider the semi-direct product
\begin{eqnarray*}
  H& := & (T\times U(1))\rtimes \Lv
\end{eqnarray*}
with multiplication
\begin{eqnarray*}
  (s,z,m) (s',z',m') & = & \(s\cdot s',z\cdot z'\cdot e^{-2\pi i I(m)}(s'),m+m'\).
\end{eqnarray*}
\begin{Cor}\label{lem:t/H}
  We have an equivalence of groupoids
  \begin{eqnarray*}
    \mT/\negmedspace/\mT^{conj} & \simeq & \ft/\negmedspace/ H
  \end{eqnarray*}
  where $T\times U(1)$ acts trivially, and $\Lv$ acts by addition. 
\end{Cor}
So, vector bundles on the inertia groupoid of $\mT$ may be thought of as $H$-equivariant bundles on $\ft$.
A source of such bundles are representations of $H$.
Writing $\CC_{\lambda,k}$
for the irreducible representation of the torus $T\times U(1)$ with
weight $(\lambda,k)$, any representation of $H$ decomposes as
\begin{eqnarray*}
  V & \cong & \bigoplus_{k\in\bbZ}V_k\tensor \CC_k,
\end{eqnarray*}
such that $U(1)$ acts on the $k$th summand with winding number $k$, and further
\begin{eqnarray*}
  V_k & \cong & \bigoplus_{\lambda\in\Lambda} V_{\lambda,k} \tensor \CC_{\lambda,k},
\end{eqnarray*}
where $T$ acts on the $\lambda$ summand with weight $\lambda$.
The action of $\Lv$ on $V$ defines isomorphisms
\begin{eqnarray*}
  \psi_m\negmedspace : V_{\lambda,k} & \xrightarrow{\,\,\,\cong\,\,\,} & V_{\lambda - kI^\sharp(m),k}.
\end{eqnarray*}
It follows that the character of $V_k$ is a linear
combination of the expressions
  $$   \sum_{m\in\Lv} e^{2\pi i \lambda-kI^\sharp(m)}.$$

\subsection{The Looijenga line bundle}
Fix a complex number $\tau$ with imaginary part $im(\tau)>0$, and set
$q=e^{2\pi i\tau}$.
Over 
\begin{eqnarray*}
  T_\CC\,/\, q^{\Lv} & \cong & \ft_\CC\,/\, (\tau \Lv + \Lv), 
\end{eqnarray*}
we have the {\em Looijenga line bundle for $I$}
\begin{eqnarray*}
  \mL_{Lo}(I) & := & (T_\CC\times \CC) \,/\, \sim ,
\end{eqnarray*}
with
\begin{eqnarray}\label{eq:sim_Looijenga}
  (h,c) \sim \(hq^m,c\,e^{-2\pi iI^\sharp(m)}(h)\cdot q^{-\frac12
    I(m,m)}\),\quad\quad m\in\Lv.
\end{eqnarray}
\begin{Lem}
  Let $\LE^I$ be the (topological) complex line bundle on $T\times T$ 
  associated to the principal $U(1)$-bundle described in Theorem
  \ref{thm:inertia_groupoid} (ii).
  Then we have an isomorphism of line bundles 
  \begin{eqnarray*}
    \LE^I & \longrightarrow & \mL_{Lo}(I).
  \end{eqnarray*}
\end{Lem}
\begin{proof}
  We use the map 
  \begin{eqnarray*}
    \ft\times\ft & \longrightarrow & \ft_\CC\\
    (x,y) & \longrightarrow & \tau x + y
  \end{eqnarray*}
  to identify the base spaces.
  To identify the line bundles, we note that the map
  \begin{eqnarray*}
    \ft \times T\times \CC & \longrightarrow &
    T_\CC\times \CC \\
    (x,t,z) & \longmapsto & \(tq^x , zq^{-\frac12I(x,x)}\)
  \end{eqnarray*}
  covers this identification and sends the equivalence relation of
  Theorem \ref{thm:inertia_groupoid} (ii) to \eqref{eq:sim_Looijenga}.
\end{proof}
\subsection{Discontinuity of the categorical character}
One could hope that the categorical character of a representation of $\mT$ was a vector bundle over $\mT/\negmedspace/\mT^{conj}$. 
In the case of the basic representation $\varrho_{bas}$ as in
Section \ref{sec:basic_representation}, however, it turns out that the
categorical character does not have the desired format: the vector space $X_{\varrho_{bas}}(t)$ does not vary continuously with $t\in T$.
For $t^*$ acting on $\mC oh(T_\bbC)$, one has
\[
  \ttr(t^*) \ =\ Nat(\id,t^*)\ = \
  \begin{cases}
    \bbC[\widehat T] & t=1\\
    0& \text{else.}
  \end{cases}
\]
Other variations of the theory exhibit the same discontinuity behaviour:
if $f\negmedspace : X\longrightarrow X$ is an automorphism of an algebraic variety over $k$, then besides the action of $f^*$ on $\mC oh(X)$ one can consider
the action on the bounded derived category $\mD^b_{coh}(X)$ or the action
by convolution with the Fourier-Mukai kernel $\mO_{\Gamma_{\negmedspace f}^t}$ associated to the (transpose) graph of $f$. It is difficult to get a handle on traces of the derived functors $Lf^*$, but
the categorical trace $\ttr(\mO_{\Gamma_{\negmedspace f}^t})$ was computed in \cite{Ganter15} and found to depend on the fixed
point locus $X^f$. This fixed point locus is empty when $f$ is given by multiplication with an element $t\in T_\bbC\setminus \{1\}$.
On the other hand, for $t=1$,
we have the trace of the identity kernel. This is given by the Hochschild cohomology group  $$\ttr^\bullet(\mO_\Delta) = HH^\bullet(T_\bbC).$$ 
Changing the definition of categorical trace to the second notion ({\em``round trace''} in \cite{Willerton2Traces}), one arrives at a similar picture with Hochschild homology in place of Hochschild cohomology.
In the setting of Section \ref{sec:connections}, where
$t^*$ acts on vector bundles with flat connection on the analytic variety $T_\bbC$, parallel transport along any piecewise smooth path from $1$ to $t$ defines an
element of $\ttr(t^*)$. Homotopic paths yield the same element.
Writing $\mL_tT_\bbC$ for the space of piecewise smooth paths from $1$ to $t$ inside $T_\bbC$, we therefore have an inclusion
\[
  \bbC\, \pi_0\(\mL_tT_\bbC\)\subseteq \ttr(t^*),
\]
and further a non-canonical bijection 
\[\Lv \cong \pi_0\(\mL_tT\).\]
In this setting, the discontinuity of $R_{bas}$ expresses itself in the non-trivial monodromy of the connection on the bundle
$\bbL_t$.

On the "twisted sector", i.e., when $\epsilon=-1$, the fixed points under $s\mapsto ts\inv$ are exactly the square roots of $t$ in $T_\bbC$.

It would be interesting to have an interpretation of our basic representation that does not suffer from such discontinuity issues.
In this context it should be noted that the homotopy theoretic counterpart to the derived fixed points $\Gamma_g^t\cap^R\Delta$ in the computation of categorical traces in the Fourier-Mukai setting are the twisted sectors of the loop space, $\mL_g X$. 
\section{Automorphisms of categorical tori}
\subsection{Automorphisms of categorical tori}
We will write $\Gamma_2\Lambda=(\Lambda\tensor\Lambda)^{S_2}$ for the second divided power of the weight lattice $\Lambda$.
This is the dual of the second symmetric power $Sym^2(\Lv)=(\Lv\tensor\Lv)/S_2$ of the coweight lattice; so elements of $\Gamma_2\Lambda$ may be thought of as symmetric bilinear forms on $\Lv$. Associated to such a symmetric bilinear form $B$ is the quadratic form 
$\phi_B(m)=B(m,m)$.
There is an exact sequence
\[
  \begin{tikzcd}
    0\ar[r] & (\Lambda\tensor\Lambda)^{sgn}\ar[r]& \Lambda\tensor\Lambda\ar[r]&\Gamma_2\Lambda\ar[r] & \Lambda/2\Lambda\ar[r] & 0,
  \end{tikzcd}
\]
where the second map is the symmetrisation map sending a bilinear form $J$ to $J+J^t$ and the third map sends a symmetric bilinear form $B$ to the class $[\phi_B]$ modulo $2$. Its kernel (the image of the symmetrisation map) consists of the even symmetric bilinear forms on $\Lv$. We will denote it by $(\Gamma_2\Lambda)_{ev}$.
The elements of $(\Lambda\tensor\Lambda)^{sgn}$ are the skew symmetric bilinear forms on $\Lv$.
Fix a bilinear form $J$ on $\Lv$ and let $I=J+J^t$ be its symmetrisation. Then an automorphism $f$ of $\Lv$ is in the isometry group $O(\Lv,I)$ if an only if $J-f^*J$ is skew symmetric. Inside the semi-direct product $O(\Lv,I)\ltimes (\Lambda\tensor\Lambda)$, consider the subgroup with elements
\[
  \widetilde{O(\Lv,I)} = \{(f,B)\mid B-B^t=J-f^*J\}.
\]
This is a non-split extension of the form
\[
  \begin{tikzcd}
    0\ar[r] & \Gamma_2\Lambda\ar[r] & \widetilde{O(\Lv,I)}\ar[r] & {O(\Lv,I)}\ar[r] & 1.
  \end{tikzcd}
\]
\begin{Def}
  We will write $E$ for the quotient group
  \[
    E\,:=\,\widetilde{O(\Lv,I)}\,/\,(\Gamma_2\Lambda)_{ev}
  \]
  and $\Xi$ for the crossed module
  \begin{eqnarray*}
    \xi\negmedspace : \Lambda & \longrightarrow & E\\
    \lambda & \longmapsto & (\id, [\lambda\tensor\lambda]), 
  \end{eqnarray*}
  where $(f,B)$ acts on $\Lambda$ via $f^*$.
\end{Def}
The sits inside a crossed module extension
\[
  \begin{tikzcd}
    0\ar[r] & 2\Lambda \ar[r]& \Lambda\ar[r,"\xi"] & E\ar[r] & {O(\Lv,I)}\ar[r] & 1.
  \end{tikzcd}
\]
Let $\mT$ be the categorical torus associated to $J$, as in Section \ref{sec:categorical_tori}.
\begin{Thm}
  \label{thm:Aut(T)}
  The categorical group of centre preserving weak automorphism of $\mT$ is equivalent to the strict categorical group $\mS(\Xi)$.
\end{Thm}
\begin{proof}
  We have an equivalence of crossed module extensions
\[
  \begin{tikzcd}
    0\ar[r] & 2\Lambda \ar[r]\ar[d,equals] &\Lambda\times_{\Lambda/2\Lambda} \Gamma_2\Lambda\ar[r,"\widetilde\xi"]\ar[d] & \widetilde{O(\Lv,I)}\ar[r]\ar[d] & {O(\Lv,I)}\ar[r]\ar[d,equals] & 1\\
    0\ar[r] & 2\Lambda \ar[r]& \Lambda\ar[r,"\xi"] & E\ar[r] & {O(\Lv,I)}\ar[r] & 1,
  \end{tikzcd}
\]
where $\widetilde\xi(\lambda,B) = (\id,B)$, and the action on the top is via the action of $O(\Lv,I)$ on $\Lambda$. Write $\widetilde\Xi$ for the crossed module on the top. The promised equivalence is given by the strict crossed module homomorphism
\begin{eqnarray*}
  A\negmedspace : \widetilde{\Xi} & \longrightarrow & wAct^+(\Theta),
\end{eqnarray*}
where $A_0$ sends the element $(f,B)$ of $\widetilde{O(\Lv,I)}$ to
the weak automorphism  
  \begin{eqnarray*}
    f_0(x) &\,=\,& f(x) \\
    f_1(m,[a]) & = & (f(m),[a]) \\
    \beta_{x,y} & = & (0,[B(x,y)]),
  \end{eqnarray*}
and $A_1$ sends the element $(\lambda,B)$ of
$\Lambda\times_{\Lambda/2\Lambda}\Gamma_2\Lambda$ to 
\[
  \eta(x) = \frac12\(B(x,x)-\lambda(x)\).
\]
To see that $A$ is an isomorphism on $\pi_1$ and hence fully faithful, 
note that
\[
  \pi_1wAct(\Theta) \,=\,\{\map\eta\ft\bbR/\bbZ\mid\eta\at{\Lv} = 0\text{ and } \eta(x)\eta(y) = \eta(xy)\}
\]
is canonically identified with $\widehat T$. Under this identification, the restriction of $A$ to $\pi_1(\widetilde{\Xi})$ becomes the standard isomorphism $\Lambda\cong\widehat T$. Since $$H^2_{gp}(T;U(1))\cong H^3(BT;\bbZ)=0,$$
we obtain from 
\eqref{eq:Z->G->Inn} an isomorphism
$$
  \pi_0(wAct^+(\Theta)) \cong O(\Lv,I).
$$ 
More precisely, the class of $A_0(f,B)$ is mapped to $f$ under the isomorphism \eqref{eq:Z->G->Inn}.
Hence $A$ is an equivalence of crossed modules.
\end{proof}
\subsection{The unimodular case}
If $I^\sharp\negmedspace:\Lv\longrightarrow \Lambda$ is an isomorphism, then the centre of $\widetilde{T[2]}$ consists of $\{\pm1\}$, and we have an isomorphism
of extensions
\[
  \begin{tikzcd}
    1\ar[r] & Inn(\widetilde{T[2]})\ar[r]\ar[d,swap,"I^\sharp_{\bbF_2}"]\ar[d,"\cong"]& Aut(\widetilde{T[2]})\ar[r]\ar[d,"\cong"] & Out(\widetilde{T[2]})\ar[r]\ar[d,"\cong"] & 1\\[+2ex]
    1\ar[r] & T^*[2]\ar[r]& E'\ar[r] & O(T[2],[\phi_I])\ar[r] & 1.
      \end{tikzcd}
 \]
 where \[E'\subset O(T[2],[\phi_I])\ltimes C^1(T[2],\bbF_2)\] is the subgroup  
\[
  E'=\{(f,c)\mid dc = {J+f^*J}\}.
\]
In this situation, $E$ is isomorphic to the pull-back of the extension $E'$
along the composite \[O(\Lv,I)\longrightarrow PO(\Lv,I)\longrightarrow O(T[2],[\phi_I]).\]
\label{eq:unimodularE}
In the case of the Leech lattice, it follows that $E$ is equal to the extension $\widetilde C$ defined on Page \pageref{page:Ctilde}.

\subsection{Automorphisms of $\mTT$}
Write $\widetilde\Lv$ for the central extension of $\Lv$ by $\bbZ/2\bbZ$ with cocycle $J^t$ (modulo $2$). It will be convenient to use additive notation for the elements of its centre. Throughout this section, $I^\sharp$ is assumned to be an isomorphism. 
\begin{Def}  
  We will write $\Xi'$ for the crossed module
    \begin{eqnarray*}
      \xi'\negmedspace :\widetilde{\Lv}/2\Lv & \longrightarrow & T\rtimes E\\
      ([n],\iota) & \longmapsto & \(\exp\(\frac n2\), \id,B\), 
    \end{eqnarray*}
    where $I^\sharp(n)\equiv\phi_B\mod 2$.
    The action of
    $T\rtimes E$ on $\widetilde\Lv/2\Lv$ is given by
    \[\([n],\iota\)^{(t,f,B)} = \([f\inv n], B(f\inv n,f\inv n)\).\]
\end{Def}

\begin{Thm}
  The categorical group of centre preserving weak automorphism of $\mT\rtimes\{\pm1\}$ is equivalent to the strict categorical group $\mS(\Xi')$. 
  \end{Thm}
  \begin{proof}
  We will write $\widetilde \Xi'$ for the crossed module
    \begin{eqnarray*}
      \widetilde\xi'\negmedspace :\widetilde\Lv\times_{\Lambda/2\Lambda}\Gamma_2\Lambda & \longrightarrow & \ft\rtimes \widetilde{O(\Lv,I)}\\
      (n,\iota,B) & \longmapsto & \(\frac n2, \id,B\), 
    \end{eqnarray*}
    where the fibred product in the source is defined using the map $I^\sharp\negmedspace : \Lv\longrightarrow \Lambda$ to map from $\widetilde \Lv$ to
    $\Lambda/2\Lambda$ and the map $\Gamma_2\Lambda\longrightarrow\Lambda/2\Lambda$ sending $B$ to $\phi_B$ modulo $2$. The action of
    $\ft\rtimes \widetilde{O(\Lv,I)}$ on $\widetilde\Lv\times_{\Lambda/2\Lambda}\Gamma_2\Lambda$ is given by
    \[\(n,\iota,S\)^{(a,f,B)} = \(f\inv n, B(f\inv n,f\inv n),f^*S\).\]
    The strict homomorphism
    \[
      \begin{tikzcd}
        \widetilde\Lv\times_{\Lambda/2\Lambda}\Gamma_2\Lambda \ar[d,swap,"\widetilde\xi'"]\ar[r]&[+2em]\widetilde\Lv/2\Lv\ar[d,"\xi'"]\\ \ft\rtimes \widetilde{O(\Lv,I)}\ar[r] &
        T\rtimes E
      \end{tikzcd}
    \]
    obtained by modding out $2\Lv\times(\Gamma_2\Lambda)_{ev}$ in degree $1$ and $\Lv\times(\Gamma_2\Lambda)_{ev}$ in degree $0$ induces isomorphisms on $\pi_0$ and
    $\pi_1$. We will now define a weak homomorphism $$(A',\kappa) \negmedspace : \widetilde\Xi'\longrightarrow wAct^+(\Theta').$$
       Given $a\in\ft$, we will write $c_a$ for the inner automorphism associated to $(a,1)$ as in Definition \ref{def:Norrie_actor}. Explicitly, $c_a$ is the strict
automorphism of $\mT\rtimes\{\pm1\}$ given by
\begin{eqnarray*}
  (c_a)_0\negmedspace : (x,\epsilon) & \mapsto& (x+(1-\epsilon)a,\epsilon)\\
  (c_a)_1\negmedspace : (m,[r]) & \mapsto& (m,[r-J(m,a)]). 
\end{eqnarray*}
Given $(f,B)\in\widetilde{O(\Lv,I)}$, we define the weak automorphism $(f,\beta)$ of $\mT\rtimes\{\pm1\}$, given by 
  \begin{eqnarray*}
    f_0: (x,\epsilon)&\mapsto& (f(x),\epsilon),
    \\f_1 \negmedspace: {(m,[r])}&\mapsto& {(f(m),[r])}\\ \beta_{(x,\epsilon),(y,\delta)}&=&{(0,[\epsilon B(x,y)])}.
  \end{eqnarray*}
  This extends the weak automorphism of $\mT$ in the proof of Theorem \ref{thm:Aut(T)}.
    We define $(A',\kappa)$ to be given by the pointed maps
    \[
      \begin{tikzcd}
        (n,\iota,B) \ar[r,mapsto] &[+1cm]\eta=\eta_{(n,\iota,B)}\\[-1.3em]
        \widetilde\Lv\times_{\Lambda/2\Lambda}\Gamma_2\Lambda \ar[r,"A_1'"]\ar[d,swap,"{\widetilde\xi'}"] & Maps_*(\ft\rtimes\{\pm1\},\Lv\times\bbR/\bbZ)^\times\ar[d,"\delta"]\\[+1.5em]
        \ft\rtimes \widetilde{O(\Lv,I)}\ar[r,"A_0'"]& wAct^+(\Theta')_0\\[-1.3em]
        (a,f,B)\ar[r,mapsto] & c_a\circ (f,\beta),
      \end{tikzcd}
    \]
    with
    \[
      \eta(x,\epsilon) = \(\frac{\epsilon-1}2n,\left[\frac12\(B(x,x)+J(n,x)+\frac{\epsilon-1}2\iota\)\right]\),
    \]
    along with the K\"unneth
    \[
      \kappa_{(a,f,B),(b,g,B')}(x,\epsilon)\,=\, (0,[\epsilon B(g(x),b) - B(b,g(x)) + (\epsilon-1)B(b,b)]).
    \]
    By a lengthy yet straight-forward computation, one checks that the pair $(A',\kappa)$ satisfies Axioms (W1) -- (W5). The Bernstein centre
    of $\widetilde \Xi'$ is given by the central subgroup of $\widetilde\Lv$, so $\pi_1(\widetilde\Xi')\cong\bbZ/2\bbZ$, while
    $$\pi_1(wAct^+(\Theta')=\Hom(\ft\rtimes\{\pm1\},\Lv\times \bbR/\bbZ) = Hom(\{\pm1\},\bbR/\bbZ).$$
    Explicitly, the map between these two groups equals
    $$A_1'(0,\iota,0)(x,\epsilon) = \(0,\left[\frac{(\epsilon-1)\iota}4\right]\)$$
    (or in multiplicative notation, $(x,\epsilon)\mapsto \epsilon^\iota$).
    So, $(A',\kappa)$ induces an isomorphism on $\pi_1$.
    To determine the effect of $(A',\kappa)$ on $\pi_0$, we use the short exact sequence
    \[
      \begin{tikzcd}
        0\ar[r]& H_{gp}^2(T\rtimes\{\pm1\},U(1))\ar[r] & \pi_0(\mA ut^+(\mT\rtimes\{\pm1\}))\ar[r] & {Stab([\alpha])}\ar[r]& 1
      \end{tikzcd}
    \]
    from \eqref{eq:Z->G->Inn}. 
  The Lyndon-Hochschild-Serre spectral sequence 
  \[
    H^p_{gp}\(\{\pm1\};H^q_{gp}(T;U(1))\) \,\, \Longrightarrow \,\,
    H^{p+q}_{gp}(\TT;U(1))
  \]
  yields the isomorphism
  \[
    H^2_{gp}(\TT;U(1)) \,\,\cong\,\,
    H^1_{gp}(\{\pm1\};\Lambda_-)\,\,\cong\,\, \Lambda/2\Lambda.
  \]
  These elements
  are the centre preserving automorphisms covering the identity of $T\rtimes\{\pm1\}$, and they are represented by the strict automorphisms
  \begin{eqnarray*}
   h_\lambda\negmedspace: (x,\epsilon) & \mapsto& (x,\epsilon)\\
    {(m,[r])} & \mapsto & {\(m,\left[r+\frac12\lambda(m)\right]\)}.
  \end{eqnarray*}
  Given $\lambda\in \Lambda$, let $B\in\Gamma_2\Lambda$ be such that $\phi_B\equiv\lambda\mod 2$, and set $\eta(x,\epsilon)= (0,[\frac12 B(x,x)])$.
  Then $\delta(\eta) = (h_\lambda,\beta)$, with $$\beta_{(x,\epsilon),(y,\delta)} = (0,[\epsilon B(x,y)]).$$ It follows that
  $A_0(0,\id,B)$ represents the same class as $h_\lambda$ in $\pi_0(wAct^+(\Theta'))$.
  Hence $\pi_0(A)$ maps $\Gamma_2\Lambda/(\Gamma_2\Lambda)_{ev}$ isomorphically
  onto $H^2_{gp}(\TT;U(1))$, and we have a map of short exact
  sequences
    \[
      \begin{tikzcd}
        0\ar[r] & \Gamma_2\Lambda/(\Gamma_2\Lambda)_{ev}\ar[r]\ar[d,"\cong"] & \pi_0(\widetilde\Xi')\ar[r] \ar[d,"\pi_0(A')"]& T\rtimes O(I,\Lv)\ar[r]\ar[d] & 1\\
        0\ar[r] & \Lambda/2\Lambda\ar[r] & \pi_0(\mA ut^+(\TT))\ar[r] & {Stab([\alpha])}\ar[r]&0.
      \end{tikzcd}
    \]
    Here $$Stab([\alpha]) = T\rtimes O(\Lv,I)\subseteq T\rtimes GL(\Lv) \cong Aut(\TT).$$
  The left-most vertical arrow in \eqref{eq:Z->G->Inn}
  is the restriction of the map $\pi_0(\mA d)$ to the centre, and it is given
  by  
  \[
    \begin{tikzcd}
      T[2] \ar[r] &
      \Lambda/2\Lambda \\[-4ex]
      \exp(\frac n2) \ar[r,mapsto] & {[I^\sharp(n)]}.
    \end{tikzcd}
  \]
  In other words, if $n$ is such that $I^\sharp(n) = \lambda$, then
  $c_{\frac{n}2}$ also represents the same class as $h_\lambda$, and the subgroup
  $T \rtimes\{\pm1\}$ of 
  $T\rtimes O(\Lv,I)$ 
  is identified, via \eqref{eq:Z->G->Inn} with $Inn(T\rtimes\{\pm1\})$.
The class of $[(a,f,B)]\in\pi_0(\widetilde\Xi')$ maps to the element $(\exp(2a),f)$ of $T\rtimes O(\Lv,I)$. It follows that the induced map on cokernels is an isomorphism, and hence, by the Five Lemma, so is $\pi_0(A')$. 
\end{proof}

\begin{Cor}
  We have
  \[
    \pi_0(\mO ut(\mTT))\,\cong\, PO(\Lv,I).
  \]
\end{Cor}
\begin{proof}
  Since $I^\sharp$ is an isomorphism, the proof of the Theorem implies that the leftmost vertical map of \eqref{eq:Z->G->Inn} is also an isomorphism, so we have an isomorphism between the cokernels of the middle and the right vertical arrows.
\end{proof}
\begin{Cor}
  The exact sequence \eqref{eq:seven-terms} for $\mTT$ takes the form
\[
    \begin{tikzcd}
      1\ar[r]& U(1)\ar[r,"\cong"]&
      U(1) \ar[r]& 1
    \end{tikzcd}
  \]
  \[
      \begin{tikzcd}
     1\ar[r]&       
      \{\pm1\}  
      \arrow[r,"\cong"]
      &\pi_0\mZ(\mTT)\ar[r]& 1
     \end{tikzcd}
   \]\[
     \begin{tikzcd}
      1\ar[r] & \TT\ar[r]&
      \pi_0\mA ut^+(\mTT)\ar[r]&
      PO(\Lv,I) \ar[r]& 1.
    \end{tikzcd}
\]
\end{Cor}
It would be interesting to see whether the cokernel of a canonical map $$\mTT\longrightarrow \mN(\varrho_{bas})$$
(with a suitable definition of
categorical normaliser $\mN(\varrho_{bas})$) can be related to the group $$C=2^{1+24}.\,Co_1$$
in the construction of the Monster.
\bibliographystyle{alpha}
\bibliography{monster.bib}

\end{document}